\documentclass[preprint,12pt]{elsarticle}

\usepackage{amssymb}
\usepackage{graphics}
\usepackage{latexsym}
\usepackage{amsfonts}
\usepackage{amsmath}
\usepackage{amsthm}
\usepackage{multicol}
\usepackage[mathscr]{eucal}
\usepackage{graphicx}
\usepackage{color}
\usepackage{epsfig}
\usepackage{array}
\usepackage{url}
\usepackage{hyperref}
\usepackage{amstext}
\usepackage{rotating}
\usepackage{subfig}

\newcommand{\be}{\begin{equation}}
\newcommand{\ee}{\end{equation}}
\newcommand{\bea}{\begin{eqnarray}}
\newcommand{\eea}{\end{eqnarray}}

\journal{Pattern Recognition}

\begin{document}

\begin{frontmatter}



\title{Locally Linear Embedding Clustering Algorithm for Natural Imagery}


\author{Lori Ziegelmeier, Michael Kirby, Chris Peterson}

\address{Department of Mathematics, Colorado State University, Fort Collins, CO 80523}

\begin{abstract}

The ability to characterize the color content of natural imagery is an important application of image processing.  The pixel by pixel coloring of images may be viewed naturally as points in color space, and the inherent structure and distribution of these points affords a quantization, through clustering, of the color information in the image.  In this paper, we present a novel topologically driven clustering algorithm that permits segmentation of the color features in a digital image.  The algorithm blends Locally Linear Embedding (LLE) and vector quantization by mapping color information to a lower dimensional space, identifying distinct color regions, and classifying pixels together based on both a proximity measure and color content. It is observed that these techniques permit a significant reduction in color resolution while maintaining the visually important features of images.

\end{abstract}


\end{frontmatter}



\section{Introduction}

Manifold learning in data analysis assumes that a set of observations, taken as a whole, is locally well approximated by a topological (or even geometric) manifold.  This assumption implies that the data is locally well approximated by a linear space, i.e., it is locally flat.  A fundamental goal of manifold learning is to uncover the underlying structure of this approximating manifold and to find low dimensional representations that preserve the structure and topology of the original data set optimally \cite{Weinberger}, \cite{Tenenbaum}, \cite{saul}.  A frequent simplifying assumption is that the local dimension is constant over the entire data set.   Alternatively, one may model a set of data as a collection of manifolds, allowing for intersections and for variations in dimension. For instance, the union of the $xy$-plane and the $z$-axis is not a manifold but decomposes naturally as a union of two manifolds of differing dimension. We have found this multiple manifold assumption to be appropriate for natural imagery consisting of distinct objects, e.g., a landscape image consisting of flowers, cacti, and ground vegetation.  

Points in a data set can typically be thought of as lying close to a low dimensional manifold if the points are parameterized by a relatively small number of continuous variables \cite{Kohonen}, \cite{Weinberger}, \cite{Donoho}.  For instance, a manifold structure could underlie a collection of images of a single object undergoing a change of state (such as illumination, pose, scale, translation, etc.). One way to uncover this structure is to map the collection to a high dimensional vector space by considering each image as a point with dimensionality corresponding to the number of pixels in the image and with coordinate values corresponding to the brightness of each pixel \cite{saul}, \cite{Huang}, \cite{Tenenbaum}. Many algorithms have been implemented on such data sets in order to uncover a low dimensional manifold that reflects the inherent structure of the high dimensional data.  Linear methods such as Principal Component Analysis \cite{Jolliffe} and Multidimensional Scaling \cite{Cox} have been around for many years while nonlinear methods such as ISOMAP \cite{Tenenbaum}, Locally Linear Embedding \cite{saul}, Hessian Eigenmaps \cite{Donoho}, and Laplacian Eigenmaps \cite{belkin} are more recent and have proven capable of extracting highly nonlinear embeddings.

In this paper, we focus on the Locally Linear Embedding (LLE) algorithm applied at the pixel level.  More precisely, our data sets do not consist of a set of images but rather the pixels comprising a single image.  Analysis of LLE applied to pixels has been implemented previously in \cite{roux} and has been considered in the context of hyperspectral images by \cite{bachmann}, \cite{Han}, \cite{chen}.  These works confirm the existence of an underlying structure. The goal of this paper is to utilize LLE to represent the underlying structure of color data in an image as a union of linear spaces and to quantize color space accordingly.
While examples will be drawn from the color space of digital images within the visible spectrum, it is important to note that such images are a special case of hyperspectral imagery. Implementations in one setting can typically be implemented in the more general setting with minor modifications.

\begin{figure}[h]
 	\centering 
 	\scalebox{0.4} 
 	{\includegraphics{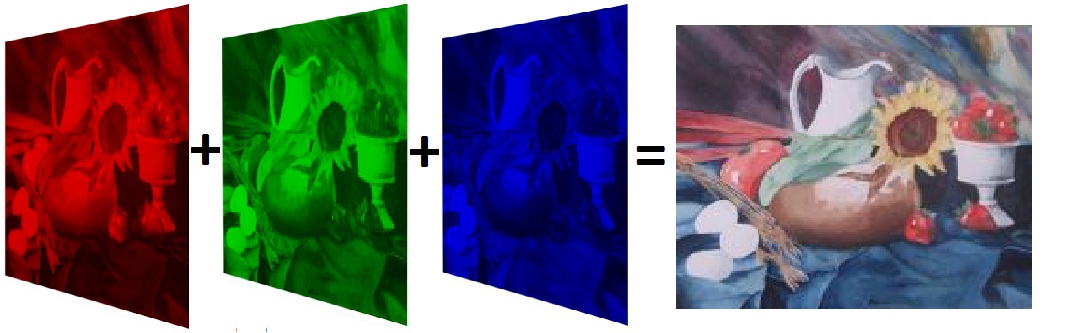}}
 	\caption{Illustration of the red, green, and blue sheets associated to pixels of an image.}
 	\label{fig:StillLifeSheets}
\end{figure}

An object under a fixed illumination condition, as perceived by the human eye, is often represented by considering a particular map to $\mathbb R^3$ obtained by integrating, at each small region of an object, the product of the spectral reflectance curve against three particular frequency response curves. We refer to these three functionals as maps to red, green, and blue (RGB) space \cite{Westland}.  As a result of this map, a digital photograph typically represents a given object/illumination pair as an $A\times B\times 3$ data array where the first two coordinates record the location in the image and the last coordinate records the values of the red, green and blue functionals on the associated spectral reflectance curve. By combining the three $A\times B$ color sheets, one can well approximate the human perception of the object/illumination pair, see Figure \ref{fig:StillLifeSheets}.

The data we will be considering consists of $A\times B\times 3$ arrays corresponding to digital pictures of natural imagery. The entries in the three $A\times B$ sheets correspond to the energy near the red, green, or blue frequencies at each pixel.  Each color component of a pixel is an integral value between 0 and 255 (corresponding to an eight bit representation).  Each pixel is associated to a point in $\mathbb{R}^3$ representing the (RGB) color of the pixel.  As there are three components of each color with 256 possible choices each, there are $256^3=16,777,216$ distinct colors that can be represented.  The color content of an image can typically be quantized at a much coarser level (while maintaining most of the visual information) by allowing each pixel's color to be identified with a prototype as determined by a quantization algorithm.  
Much work has been done to quantize the color space of natural imagery. The wide variety of approaches include statistical-based, graph theoretical, clustering, gradient descent, among many other techniques \cite{Lee}, \cite{Huang}, \cite{Cheng}, \cite{Papamarkos}, \cite{HDCheng}, \cite{Heckbert}, \cite{Orchard}, \cite{Linde}, \cite{Deng}, \cite{Velho}, \cite{ng}.  

In the LLE algorithm, a weighted graph is first constructed from a set of data as a stand-in for the local manifold structure \cite{saul}. The algorithm next determines a set of $d$ {\it embedding vectors} by discarding the eigenvector corresponding to the smallest eigenvalue of an associated graph Laplacian and keeping the $2^{nd}$ through $d+1^{st}$ eigenvectors. Arranging the $d$ vectors as columns of a matrix,  the rows of this matrix provide a map of the original data to $\mathbb R^d$. The graph Laplacian encodes the number of connected components of the graph as the dimension of its null space. Interpreting results in the LLE algorithm becomes problematic if the null space has dimension greater than one as eigenvectors in a null space are unique only up to rotation.  Thus, a canonical ordering of eigenvectors is ill defined.  As it is quite reasonable to expect the color space of natural images to be lying on multiple manifolds, it is also natural to expect multiple connected components amongst the union of the manifolds. In order to alleviate this problem, we perturb the graph Laplacian in the direction of a circulant graph Laplacian to reduce the co-rank to one.
From the $d$-dimensional embedding, we apply a technique for proximity/color segmentation.  This is done through an unsupervised clustering algorithm that exploits the topology preserving properties of LLE. 

Put another way, natural images are not random; they have structure in their color space in that adjacent pixels tend to have similar color values. These piecewise continuous variations lead to a piecewise manifold approximating the data. Through LLE, this piecewise manifold is revealed as a piecewise linear manifold. The segmentation is accomplished by uncovering the principal direction of an epsilon ball of points and segmenting the data such that points determined to be close enough to this principal vector and similarly colored to the center of the epsilon ball are classified together and removed from the data.  This iterative approach has proven robust in the presence of noise, with the input parameters reflecting the accuracy of segmentation desired.   
Thus, the algorithm exploits the transformation of local one-manifold structure to local linear structure in the mapped data.  

In Section \ref{sec:lle}, we present a brief overview of the Locally Linear Embedding algorithm and present the graph Laplacian perturbation to reduce to the case of co-rank one.  Section \ref{sec:llec} discusses the algorithm in conjunction with color quantization.  We present an example to observe that the geometric structure of a color image is revealed in a reconstruction image by exploiting locally-linear variations in pixel space through subspace segmentation.  Section \ref{sec:implementation} demonstrates the algorithm's ability to reduce the color space of natural imagery and uses this algorithm in conjunction with the classical Linde-Buzo-Gray vector quantization algorithm \cite{kirby}, \cite{Linde} in the context of a landscape ecology application. The contributions of this paper include a method for resolving LLE rank issue problems (without carrying out a decomposition into connected components) in such a manner that the local topological structure of the data is preserved, a technique for subspace segmentation, and the development of an associated clustering algorithm.

\section{Connecting Components in Locally Linear Embedding} \label{sec:lle}

\subsection{Locally Linear Embedding Algorithm}\label{sec:llealg}
The Locally Linear Embedding (LLE) algorithm \cite{saul} is an unsupervised dimensionality reduction algorithm that determines a mapping of data, lying in a higher dimensional vector space, to a lower dimensional vector space while optimizing the maintenance of local spatial relationships within the data.  Through this map, the LLE algorithm uncovers a lower dimensional representation of the data with the goal of preserving the topology and neighborhood structure of the original higher dimensional data.  The first step of the algorithm requires a criterion to determine the nearest neighbors of each data point.  The second step is to associate to each neighbor a weight.  This weight is calculated by solving a least squares problem that minimizes a certain reconstruction error $\epsilon (W)$.  More precisely, if $\textbf{x}_i$ denotes the $i^{th}$ data point from a set of $p$ points and $N_i$ denotes the indices of its set of nearest neighbors, one determines the values of $w_{ij}$ that minimize the expression
$$\epsilon(W) = \displaystyle \sum_{i=1}^p \| \textbf{\textbf{x}}_i-\displaystyle\sum_{j\in{N_i}}w_{ij}\textbf{x}_j \|^2$$
The final step of the algorithm is to determine a set of lower dimensional vectors, $\textbf{y}_i$, that minimize the function
$$\phi(\textbf{Y}) =\displaystyle \sum_{i=1}^p \| \textbf{y}_i-\displaystyle \sum_{j\in{N_i}}w_{ij}\textbf{y}_j \|^2$$ 
Let $W$ denote the $p\times p$ matrix containing the weights $w_{ij}$ (padded out with zeros). The $y_i$'s are found by solving the eigenvector problem $MY^T=Y^T\Lambda$ where $M = I-W-W^T+W^TW = \left(I-W\right)^T\left(I-W\right)$, $\Lambda$ is the diagonal matrix of Lagrange multipliers, and the $i^{th}$ row of $Y^T$ corresponds to $y_i$. 

If our data set corresponds to a sampling of a manifold and if this sampling is sufficiently dense, then a fundamental assumption of the algorithm is that each data point and its nearest neighbors can be characterized by a locally linear patch of the manifold, hence the name Locally Linear Embedding.  Data points that were close together in the original higher dimensional space should still be close together after mapped to lie in the lower dimensional space thus preserving the topology of the original data set.  Details of the LLE implementation can be found in \ref{app:lle}.  

\subsection{Example} \label{sec:lleex}
\begin{figure}[h]
 	\centering 
 	\scalebox{0.4} 
 	{\includegraphics{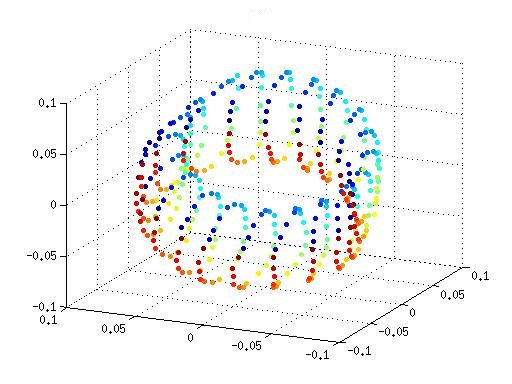}}
 	\caption{A plot of the embedding vectors obtained by LLE of images as described above. }
 	\label{fig:1sttorus}
 \end{figure}

This section consists of an example illustrating the LLE algorithm's topology preserving capabilities. The data set consists of images of a black square translated over a background of random noise. The black square is allowed to split and ``wrap'' around each boundary edge of the background.  To construct the data set, we start with a $20 \times 20$ matrix  consisting of random entries between 0 and 1.  Within this matrix, a $10 \times 10$ zero matrix is superimposed.  The data set is generated by considering all positions of the $10 \times 10$ matrix inside the $20 \times 20$ matrix of random noise (allowing both horizontal and vertical wrapping).  Thus from a topological point of view, the data set corresponds to a noisy sampling of a torus in $\mathbb{R}^{400}$.  We defined the nearest neighbors, of each element in the data set, to be the $4$ nearest data points in $\mathbb{R}^{400}$. The LLE algorithm was then used to map the data to $\mathbb R^3$.  The resulting embedded data in $\mathbb{R}^{3}$ is displayed in Figure \ref{fig:1sttorus} and reflects the original topological structure quite clearly.
 In general, results from experiments suggest that LLE can indeed be successful in its goal of a nonlinear dimensionality reduction that captures inherent topological properties through a single, linear algebra derived, map.  Due to its effectiveness, we have chosen to exploit this reduction in the context of a clustering algorithm.

\subsection{Special Considerations in Implementation of LLE on Natural Imagery}

As natural images have the feature that many pixel colors are quite similar, it would not be surprising to find pixels with identical colors.  Thus, when we consider an image as a collection of points in $\mathbb{R}^3$, we may observe distinct image pixels whose distance apart is zero. In determining nearest neighbors for a point $x_i$, we have opted to only include points whose distance from $x_i$ is greater than zero.  
We have observed that for many natural images, the $k$ nearest neighbor's graph has corank larger than 1 indicating more than one connected component.  For disconnected data,  LLE can be implemented on each of the graph's connected components separately \cite{saul}.   In this paper we have chosen a different (and slightly unusual) path in that we artificially connect components by perturbing in the direction of a cycle.  

\subsection{Connecting Disconnected Components} \label{sec:nullspace}

The Laplacian of a graph has 0 as an eigenvalue with multiplicity equal to the number of connected components of the graph \cite{Chung}.  In a similar manner, the matrix $M$, in the final step of the LLE algorithm, has co-rank corresponding to the number of connected components within the data set where connections are made by linking each data point with its nearest neighbors.  Choosing the number of nearest neighbors to be small can lead to many disconnected components.
As previously stated, while \cite{saul} indicates that LLE be implemented on each of the graph's connected components separately, we have chosen to proceed down a different path by artificially connecting previously disconnected components through a perturbation (much like the second step in LLE that adds a regularization term to the covariance matrix $C$ that would be singular in the case when $k>D$).  Here, we perturb $M$ in the direction of a matrix $T$ that has a similar structure to $M$ in that it is positive semidefinite with row sums equal to 0. We chose $T$ to be the Laplacian matrix of a cycle, i.e.
$$T=\begin{bmatrix} 1 & -1 & 0 & &\cdots & 0 \\
								-1 & 2 & -1 & 0 & \cdots & 0 \\
								 &  & \ddots &  \\
								0 & \cdots & 0 &-1 & 2 & -1 \\
								0 & \cdots & &0 & -1 & 1
\end{bmatrix}
 $$ 
Thus, in practice, we consider $M^\prime=M+\lambda T$ where $\lambda$ is a scalar.  Matlab experiments suggest that using $\lambda=10^{-9}$ produces a matrix $M^\prime$ that is artificially connected (i.e. has a corank of 1).  As $T$ comes from a cycle, the eigenvectors of $M^\prime$ inherit this circular structure while retaining much of the original topology of the pixel data.

\section{Locally Linear Embedding Clustering}\label{sec:llec}

In this section, we consider a novel way to extend the LLE algorithm to a clustering algorithm.  In natural imagery, the variation in color hues of neighboring pixels is often slight. Furthermore, pixels which are not spatially close may also exhibit very slight color variations. We can associate these slight color changes with continuous variations of pixel colors in an approximating manifold. Topological structure associated with these continuous variations is revealed by uncovering multiple underlying manifolds. These manifolds are detected using the Locally Linear Embedding algorithm.  Segmenting based on these manifolds allows for quantization of the color space, leading to a clustering algorithm.

\subsection{Clustering}
 Data clustering is the name given to creating groups of objects, or clusters, such that objects in one cluster have a shared set of features whereas objects in different clusters have less similarity with respect to these features.  Clustering is a fundamental approach to segmenting data. There are many ways to attach pairwise similarity scores to a set of data points and it is important to realize that the clusters could have very different properties and/or shapes depending on these scores.  Clustering can be done in a hierarchical manner where clusters are determined by using previously established clusters or in a partitioning manner where the data is clustered simultaneously into disjoint sets.  
In semi-supervised or constrained clustering algorithms, additional information such as data labels, information about the clusters themselves, etc. is available and utilized \cite{basu}, \cite{kirby}.  Unsupervised clustering algorithms, in which data is organized without any information like data labels, however, can identify major characteristics or patterns without any supervision from the user.  The algorithm of this paper, described in the following subsections, is non-hierarchical and unsupervised.

\subsection{Why Locally Linear Embedding?}
The Locally Linear Embedding algorithm, as discussed in Section \ref{sec:lle} and \ref{app:lle}, is a dimensionality reduction algorithm that can help uncover the inherent structure and topology of higher dimensional data by determining a map to a lower dimensional space that optimizes for neighborhood relationships. While LLE was intended for the purpose of revealing topological structure, the creators of this algorithm indicate in \cite{saul} that some of the ideas presented in LLE, namely the first and third steps, are similar to those of the Normalized Cut algorithm discussed in \cite{shi} and other clustering methods such as the one discussed by \cite{ng}.   

Through experimentation, it was observed that if a natural image, considered as a set of RGB color points, is embedded in $\mathbb{R}^2$ using nearest neighbor sets of size $4$ and a perturbation matrix $T$ with $\lambda=10^{-9}$ (as discussed in subsection \ref{sec:nullspace}), then the data lies on a relatively small collection of lines.  When the inherent color of each of the higher-dimensional input vectors was superimposed on the corresponding reconstruction vectors, it was noticed that similar colors fell along the same line.  We observed simpler (but similar) behavior when considering a fixed pixel in a set of images of a fixed object under changing ambient illumination conditions. In each of these two cases, the lines in the reduced space suggest a method for quantization. 

\begin{figure}
  \centering
  \subfloat{\label{fig:light} \includegraphics[scale=.23]{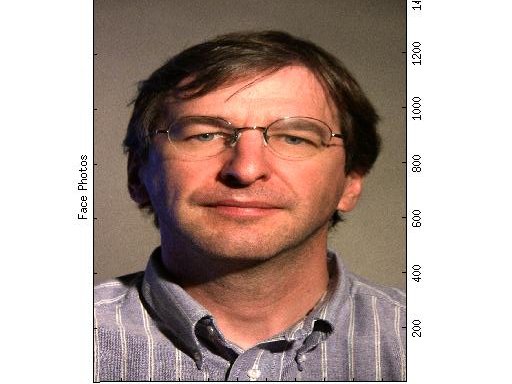}}  
  \subfloat{\label{fig:medium} \includegraphics[scale=.23]{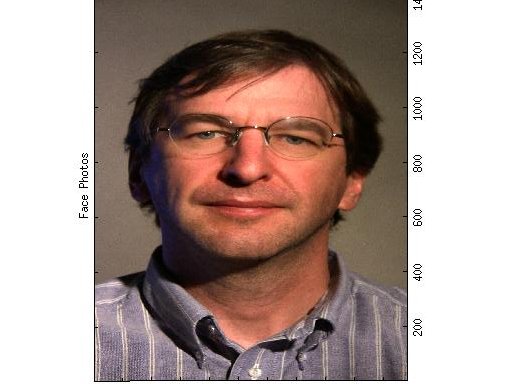}}
  \subfloat{\label{fig:dark} \includegraphics[scale=.23]{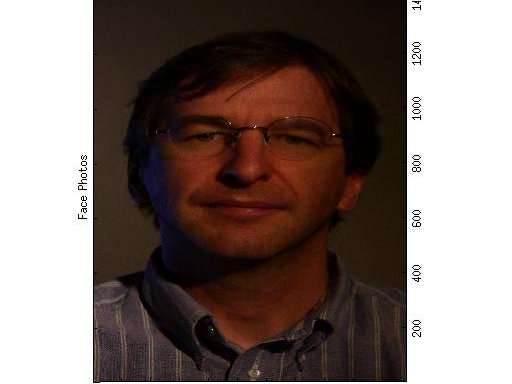}} 
  \caption{Images generated by Pattern Analysis Laboratory at Colorado State University where an individual remains motionless and the illumination of the surrounding area varies.}
  \label{fig:lighting}
\end{figure}

The following example uses a set of images from the Pattern Analysis Laboratory (PAL) database at Colorado State University. In the data, an individual remained motionless as the ambient illumination conditions were altered (with lights of fixed spectral characteristics).  Figure \ref{fig:lighting} shows three such images with different illumination conditions.
A data set was formed by considering the RGB values of a single, fixed pixel extracted from 200 such images.  

The LLE algorithm was implemented on this data set and mapped into a 2-dimensional space using $k=4$ nearest neighbors.  The null space of $M$ turned out to be 4-dimensional indicating 4 connected components.  Thus, our choice of nearest neighbor has artificially disconnected this data set that intuitively should be connected given that it is the smooth variation of illumination at a fixed point.  Therefore, either we need to increase the number of nearest neighbors or perturb the data in such a way that the data is reconnected. Note that the smallest value of $k$ that yields a 1-dimensional null space for $M$ is $k=8$.  

\begin{figure}
  \centering
  \subfloat[Data in $\mathbb{R}^3$]{\label{fig:purple3D} \includegraphics[scale=.358]{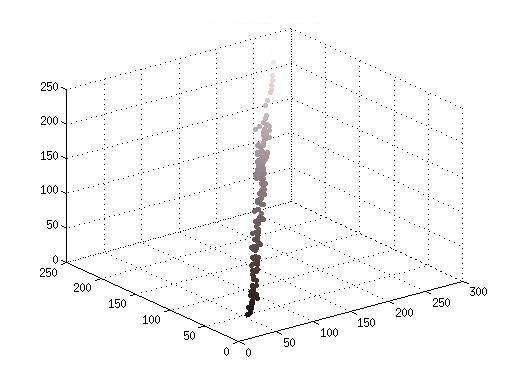}}  \hspace{-.6cm}
  \subfloat[Data in $\mathbb{R}^2$, k=8]{\label{fig:purple2D} \includegraphics[scale=.235]{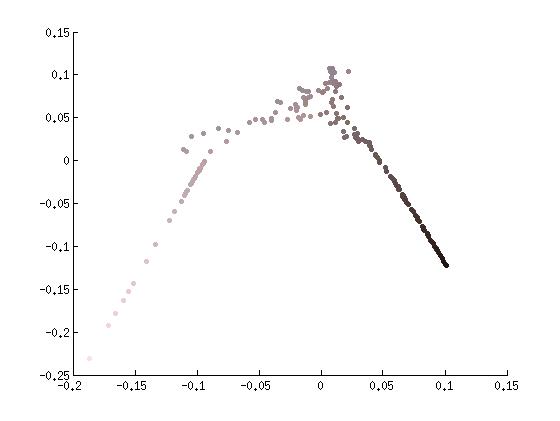}}
  \hspace{-.6cm}
  \subfloat[Data in $\mathbb{R}^2$, k=4]{\label{fig:purple2D4} 	    \includegraphics[scale=.235]{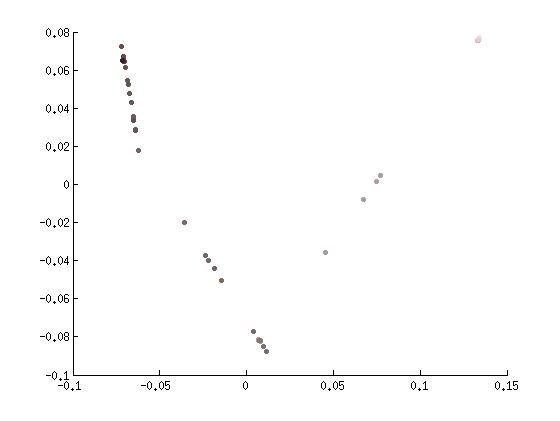}}
  \caption{Plots of 3D data points generated by extracting the RGB values of a single pixel for each of 200 images and their corresponding embedding vectors as reconstructed by the LLE algorithm using k=8 nearest neighbors, and k=4 nearest neighbors with `connected' data points.}
  \label{fig:purplepointlle}
\end{figure}

In Figure \ref{fig:purplepointlle}, we see the original plot of the RGB data points, a plot of the embedding vectors using $k=8$ nearest neighbors, and a plot of the embedding vectors using $k=4$ nearest neighbors with the perturbation discussed in Section \ref{sec:nullspace} that artificially reconnects the data set.  
We observe that the original three dimensional data appears relatively linear.  Using $k=8$ nearest neighbors, there is a degradation of this linear structure.  However, using $k=4$ nearest neighbors, with $M$ perturbed by $T$, preserves the linear structure at a local level.

Therefore, we have observed two important properties of LLE.  First, if a data set is disconnected using a choice of $k$ nearest neighbors, it can be artificially reconnected using an appropriate perturbation.  Provided the perturbation is not too extreme, this does not affect the local topology of the data as expressed by LLE, as the linear structure of similar colors falling along one dimensional subspaces holds for each component.  Second, the LLE algorithm is able to uncover the gradation of hue or variance of illumination within an image which corresponds to a linear structure in the plot of the reconstruction image of an RGB color space.  Using this linear structure, in which each data point of the reconstruction image is colored according to its corresponding high-dimensional data point, a color space clustering algorithm is obtained.  

Essentially, the Locally Linear Embedding Clustering (LLEC) algorithm segments the distinct linear manifolds that appear in the reconstruction plot of the LLE algorithm and then identifies which data points lie close to which subspace.  The RGB color information is then overlaid onto the reconstruction data and further used to segment the data by clustering similarly colored points together.  

\subsection{Subspace Segmentation}
We propose a subspace segmentation technique in the LLEC algorithm that involves selecting a point, $\textbf{y}^*$, in the reconstruction data and then constructing an epsilon ball of appropriate size centered around this point.  The point $\textbf{y}^*$ may be chosen randomly or with a more deterministic criterion such as those discussed in \ref{app:RandPoint}.  A data matrix, $A$, is formed by inserting each embedding vector falling inside this epsilon ball into the rows of the matrix.  The singular value decomposition, $A=U \Sigma V^T$, is computed to determine the right singular vector corresponding to the largest singular value.  This singular vector corresponds to the line that passes through the mean and minimizes the sum of the square Euclidean distances between the points in the epsilon ball and the line, indicating the principal direction of these data points \cite{Jolliffe}.  This unveils the 1-dimensional subspace that we are after.  Note a method for $k$-dimensional subspace segmentation is discussed below.  In order to segment the data via its proximity to each subspace, we then determine which points $\textbf{y}_i$ in the data set satisfy $||\textbf{y}_i-\mathbb{P}\textbf{y}_i||< \epsilon_1$ where $\epsilon_1$ is some tolerance to identify those points that are `close enough to' the subspace in consideration, and $\mathbb{P}$ is the projection onto the first right singular vector.   

Now, as many lines overlap or intersect, simply using proximity to the subspace will not yield an appropriate segmentation.  Thus, the data is further segmented by identifying those points, $\textbf{y}_i$, whose RGB color values, $\textbf{x}_i$, are most similar to the color of $\textbf{y}^*$, $\textbf{x}^*$, by computing $||\textbf{x}_i-\textbf{x}^*||<\epsilon_2$ where $\epsilon_2$ is again some tolerance to identify which points are `similarly' colored to $\textbf{y}^*$. Points that are identified as being `close enough to' the subspace generated by the random point and `similarly' colored to this point are identified together and removed from the data set.  The process is repeated until all points have been identified with a distinct subspace using both proximity and color information.  Thus, we obtain a clustering based on both color and spatial proximity metrics (in the 2-D representation).  

\begin{figure}
  \centering
  \subfloat{\label{fig:SubSegIll} \includegraphics[scale=.50]{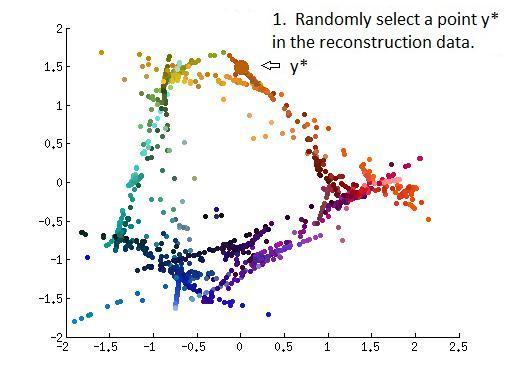}} 
  \subfloat{\label{fig:SubSegBall} \includegraphics[scale=.50]{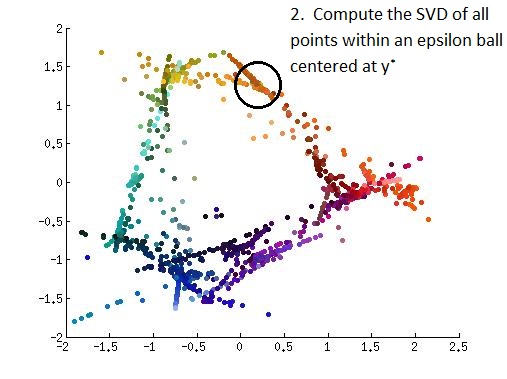}}\\
  \subfloat{\label{fig:SubSegLine} \includegraphics[scale=.50]{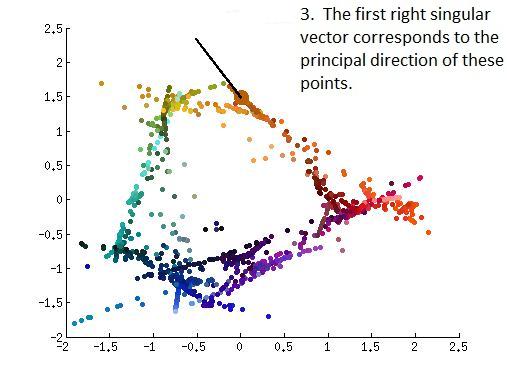}}  
  \subfloat{\label{fig:SubSegRem} \includegraphics[scale=.50]{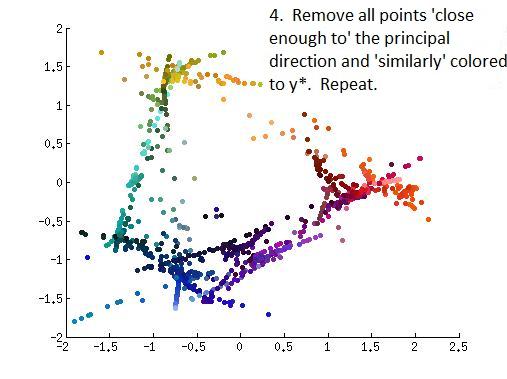}} \\ 
  \caption{Illustration of subspace segmentation algorithm.}
  \label{fig:SubSeg}
\end{figure}

\begin{figure}
  \centering
  \subfloat[Original]{\label{fig:flowers2} \includegraphics[scale=0.25]{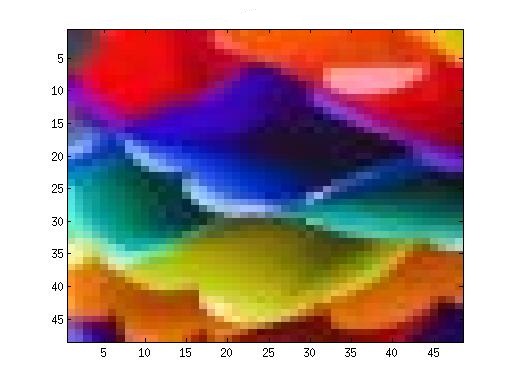}}               
  \subfloat[2D Reconstruction after LLE]{\label{fig:flowers2D2}\includegraphics[scale=0.2]{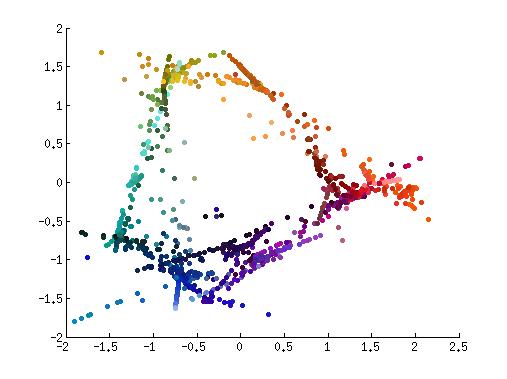}}\\
 	\subfloat{\label{fig:ss1} \includegraphics[scale=0.18]{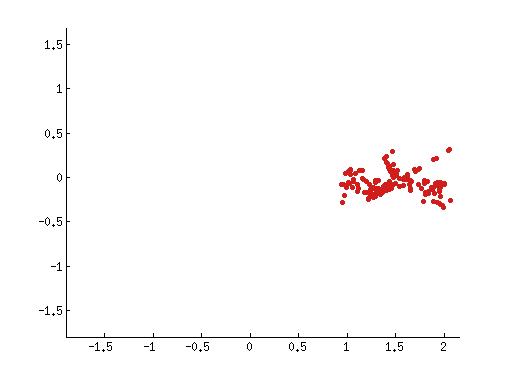}}  
 	\subfloat{\label{fig:ss2} \includegraphics[scale=0.18]{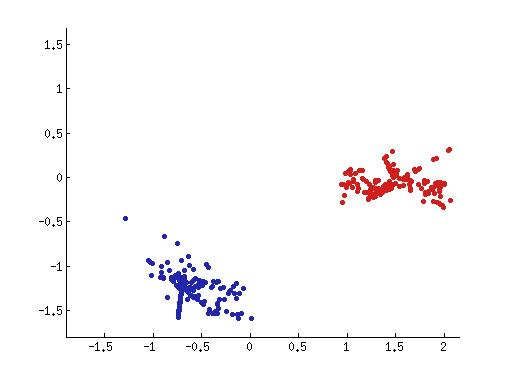}}
 	\subfloat{\label{fig:ss3} \includegraphics[scale=0.18]{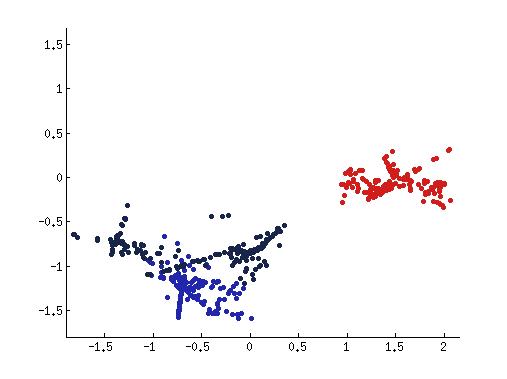}}
 	\subfloat{\label{fig:ss4} \includegraphics[scale=0.18]{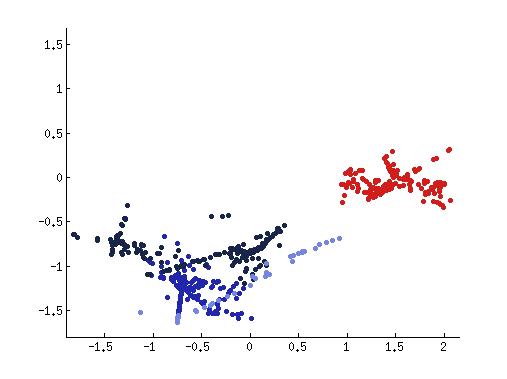}}\\
 	\subfloat{\label{fig:ss5} \includegraphics[scale=0.18]{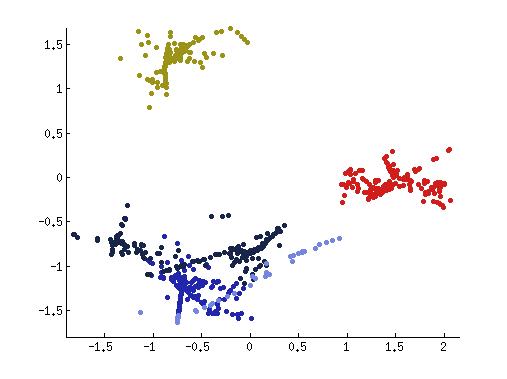}}
 	\subfloat{\label{fig:ss6} \includegraphics[scale=0.18]{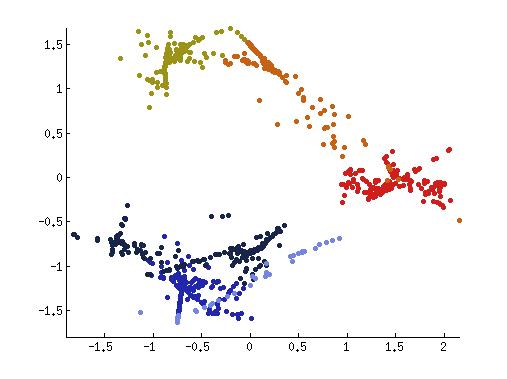}}
 	\subfloat{\label{fig:ss7} \includegraphics[scale=0.18]{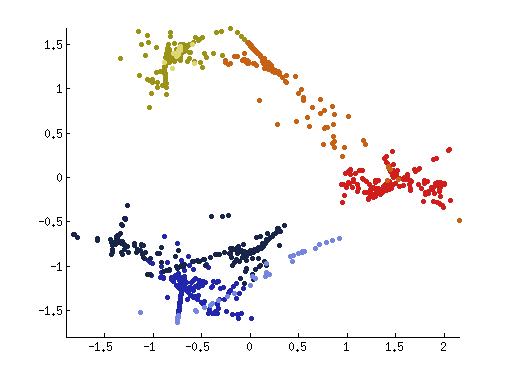}}
 	\subfloat{\label{fig:ss8} \includegraphics[scale=0.18]{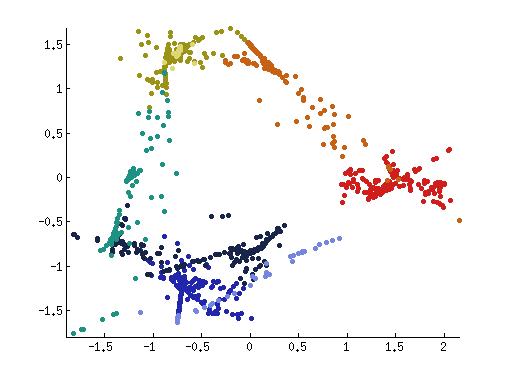}}\\
 	\subfloat{\label{fig:ss9} \includegraphics[scale=0.18]{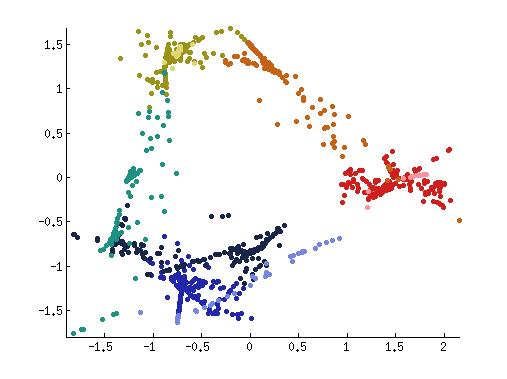}}
 	\subfloat{\label{fig:ss10} \includegraphics[scale=0.18]{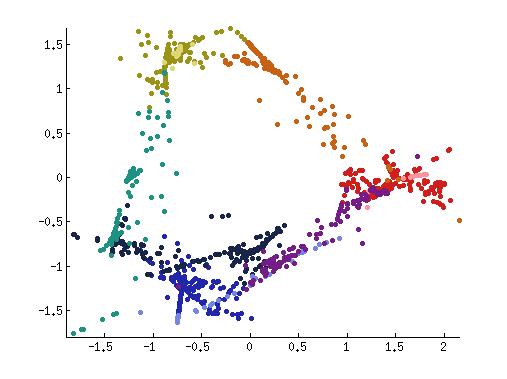}}
 	\subfloat{\label{fig:ss11} \includegraphics[scale=0.18]{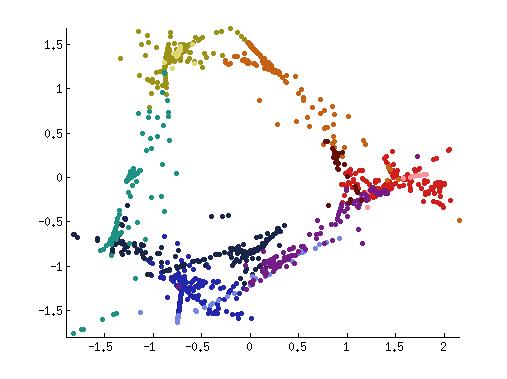}}
 	\subfloat{\label{fig:ss12} \includegraphics[scale=0.18]{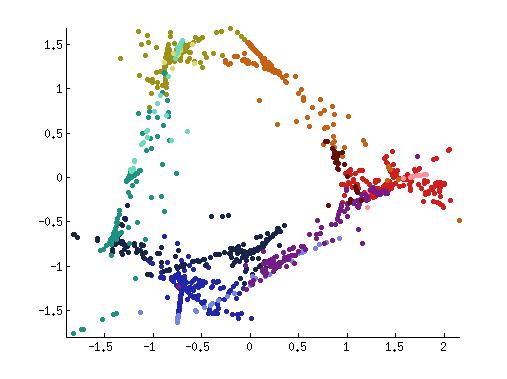}}\\
 	\subfloat{\label{fig:ss13} \includegraphics[scale=0.18]{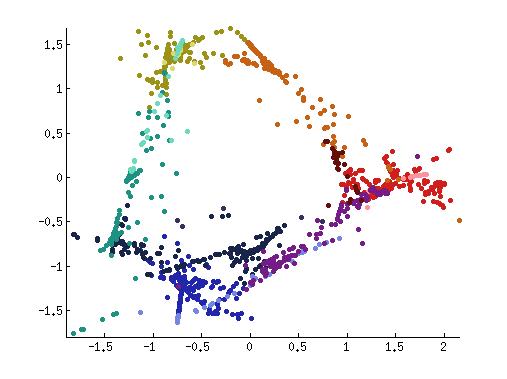}}
 	\subfloat{\label{fig:ss14} \includegraphics[scale=0.18]{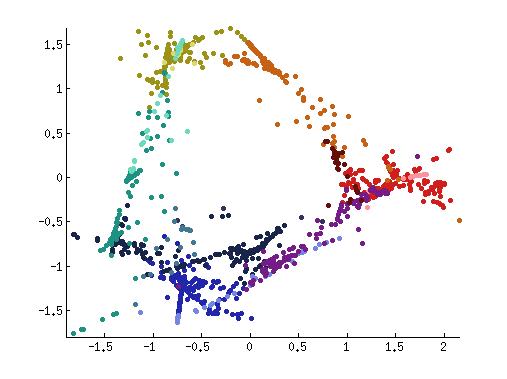}}
 	\subfloat{\label{fig:ss15} \includegraphics[scale=0.18]{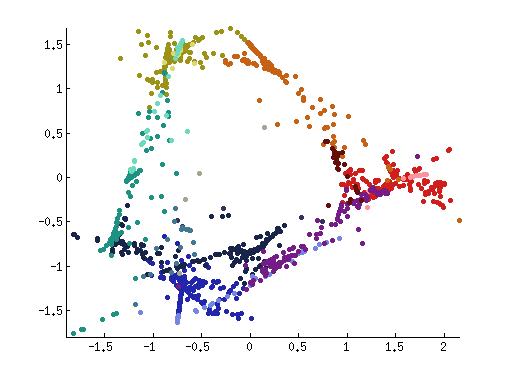}}
 	\subfloat{\label{fig:ss16} \includegraphics[scale=0.18]{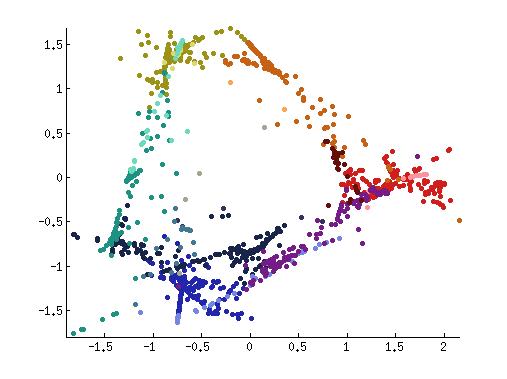}}\\
 	\subfloat{\label{fig:ss17} \includegraphics[scale=0.18]{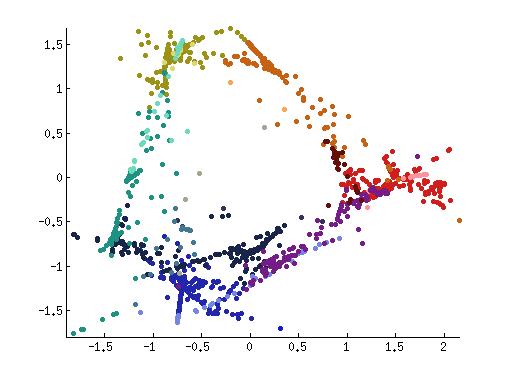}}
 	\subfloat{\label{fig:ss18} \includegraphics[scale=0.18]{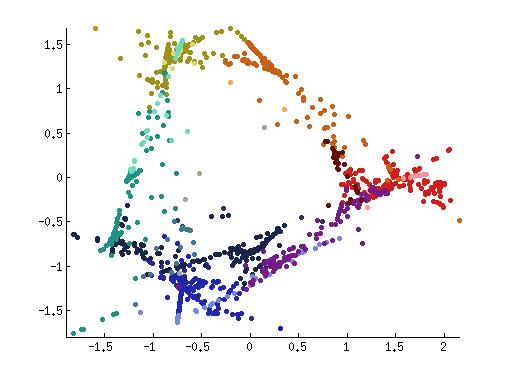}}
 	\subfloat{\label{fig:ss19} \includegraphics[scale=0.18]{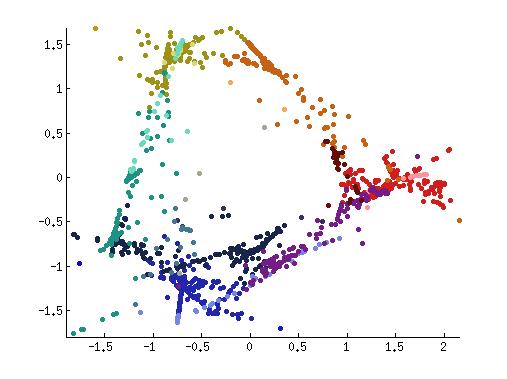}}
 	\subfloat{\label{fig:ss20} \includegraphics[scale=0.18]{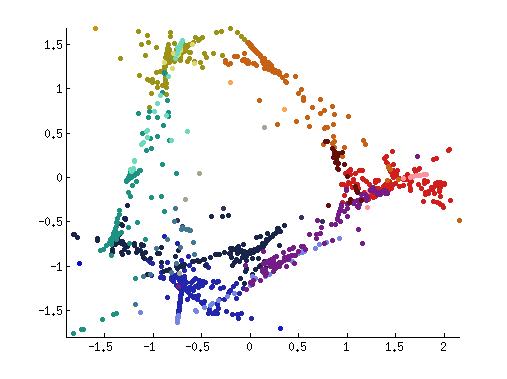}}\\
  \caption{Subspace segmentation of 2D reconstruction data using LLEC with accuracy tolerances of approximately 0.4 by initializing $\textbf{y}^*$ as the point whose $b^{th}$ nearest neighbor has the smallest distance.}
  \label{fig:subspace}
\end{figure}

Figure \ref{fig:SubSeg} provides a step-by-step illustration of this subspace segmentation algorithm.
Also refer to Figure \ref{fig:subspace} to see an actual implementation of this algorithm.  Observe that the subimages represent the subspaces iteratively removed from the reconstruction data.

We have described the approach for 1-D subspace segmentation, but this method can be used for multiple dimensions as well by considering pixels sufficiently correlated with the first several principal components.  Thus, if an $m$-dimensional subspace segmentation is desired, the $m$ right singular vectors corresponding to the $m$ largest singular values form the principal directions of the data and span the subspace we want to uncover.  Once this subspace is uncovered, the rest of the approach described above is analogous for a multidimensional subspace segmentation.  Note that this subspace segmentation technique has proven robust in the case of noise whereas methods such as Generalized Principal Component Analysis (GPCA) \cite{Vidal} have proven less effective in our experiments.

\subsection{The Locally Linear Embedding Clustering Algorithm} \label{sec:llecalgorithm}

The LLEC algorithm for quantizing the color space of natural imagery is summarized in the steps below.  First, embed a data set $X$ with points of dimension $D=3$ into a lower dimension $d=2$ using $k=4$ nearest neighbors by using LLE.  Next, identify the distinct subspaces appearing in the reconstruction data, $Y$, of the LLE algorithm by using the process described in the previous section.  Through this subspace segmentation, determine the Voronoi sets, $S_i$, formed by identifying those points that are `close enough to' each subspace and `similarly' colored to the point $\textbf{y}^*$.  Note that the number of distinct Voronoi sets is the number of distinct subspaces.  Let's call this number $S$.  Then, calculate the mean of the colors, $\mu_i = \frac{1}{|S_i|}\sum_{\textbf{y} \in S_i} \textbf{y}$ of each set.  Finally, identify all points $\textbf{y} \in S_i$ by the prototype $\mu_i$.  Note that each $\textbf{y}_i$ in the reconstruction data corresponds to a unique $\textbf{x}_i$ in the original data space, so this determines a clustering of the data set $X$.

\section{Implementation}\label{sec:implementation}

As indicated in \cite{Deng} a major complication in color space quantization often relates to varying shades of a given color due to illumination.  We have observed that the LLEC algorithm handles this illumination component by identifying various shades of a given hue as a unique subspace and all pixels that are elements of this subspace can be identified together.

At this time, the LLEC algorithm is not a fast algorithm as the procedure for performing LLE and the search to determine those points that are `close enough to' each subspace and `similarly' colored to the random point being considered are computationally intensive.  We will see however that LLEC does an excellent job of quantizing the color space of natural imagery.  A benefit of LLEC is that the only free parameters in the algorithm are $\epsilon_1$ and $\epsilon_2$, the tolerances which can be specified by the user to reflect the desired accuracy of the quantization.  The value of LLEC then is that it can be implemented on an image to determine the natural subspaces of the color space.  The knowledge obtained by segmenting these subspaces can then be used in conjunction with other clustering algorithms such as the Linde-Buzo-Gray (LBG) \cite{Linde} vector quantization algorithm to identify the starting centers as the subspaces unveiled in the LLEC algorithm.  We will see this applied shortly.

\subsection{LLEC used to Quantize Color Space}
First, let's consider the ability of the LLEC algorithm to quantize the color space of a variety of images.  In each of these examples, the images were processed using MATLAB.  Each sheet pertaining to the red, green, or blue component of pixels within an image was converted from a matrix of dimension equal to the resolution of each image to a long row vector of dimension $1 \times p$ where $p$ is the number of pixels in the image.  A new matrix, $X$, of dimension $3 \times p$ was created to contain all of the data entries of these long row vectors.  By organizing the data this way, we see that each column of the matrix corresponds to the RGB components of an individual pixel which is a data point to be analyzed.  We have chosen to use the Euclidean metric to calculate distance--a measure of proximity--between points.  

\begin{figure*}[t]
  \centering
  \subfloat[{\tiny{Original}}]{\label{fig:f1} \includegraphics[scale=.24]{flowersMatlab.jpg}}       
  \subfloat[{\tiny{$S=13$, $DE=4.2162\times10^3$}}]{\label{fig:f2} \includegraphics[scale=.183]{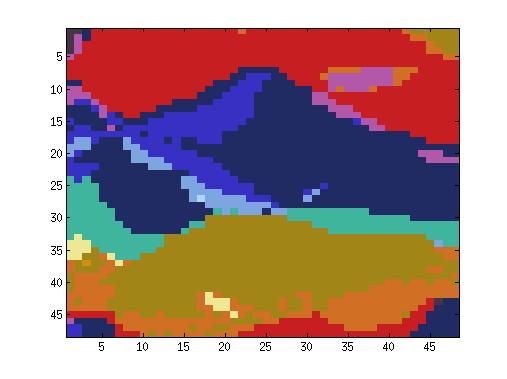}}
 \subfloat[{\tiny{$S=24$, $DE=2.2738\times10^3$}}]{\label{fig:f3} \includegraphics[scale=.183]{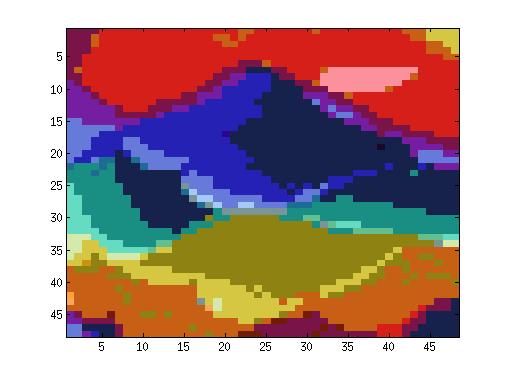}}
 \subfloat[{\tiny{$S=92$, $DE=764.723$}}]{\label{fig:f4} \includegraphics[scale=.183]{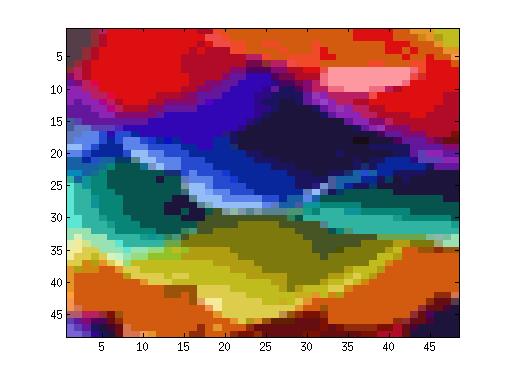}}\\
 \vspace{-.3cm}
  \subfloat[{\tiny{Original}}]{\label{fig:f5} \includegraphics[scale=.24]{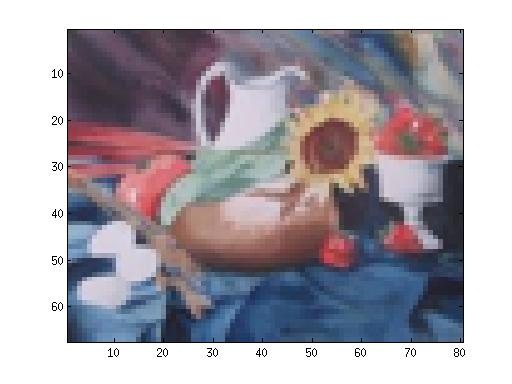}}
  \subfloat[{\tiny{$S=10$, $DE=2.0091\times10^3$}}]{\label{fig:f6} \includegraphics[scale=.183]{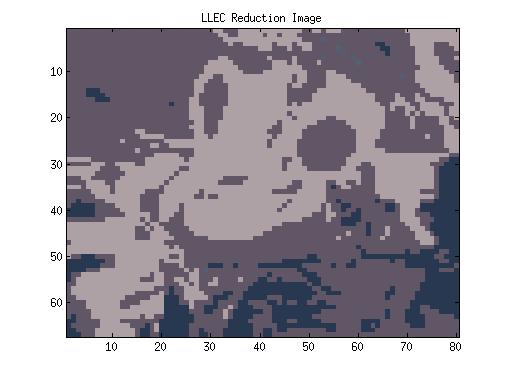}}
  \subfloat[{\tiny{$S=15$, $DE=1.4704\times10^3$}}]{\label{fig:f7} \includegraphics[scale=.183]{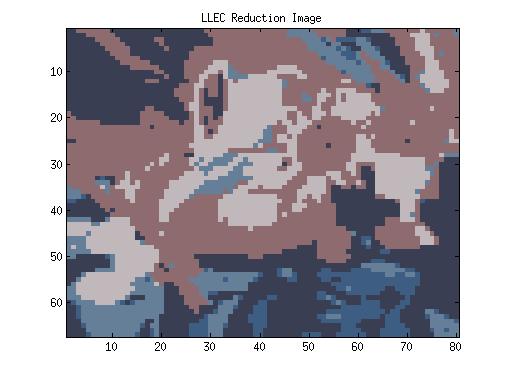}}  
  \subfloat[{\tiny{$S=57$, $DE=580.750$}}]{\label{fig:f8} \includegraphics[scale=.183]{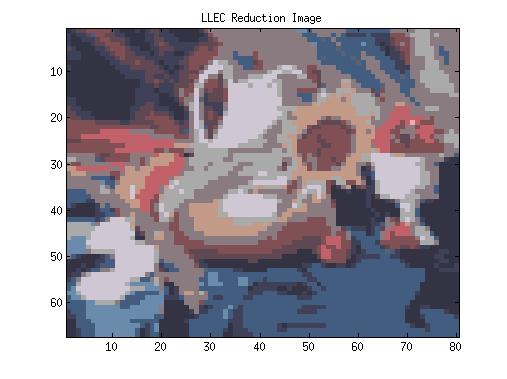}}\\
  \vspace{-.3cm}
  \subfloat[{\tiny{Original}}]{\label{fig:f9} \includegraphics[scale=.24]{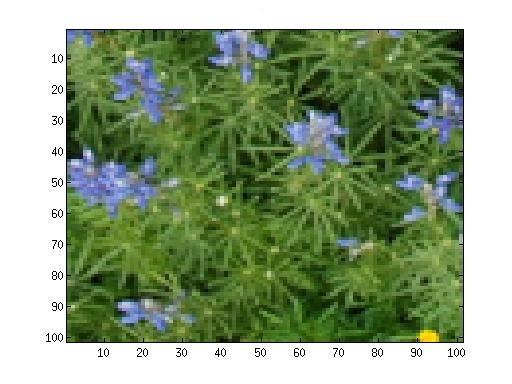}}
  \subfloat[{\tiny{$S=5$, $DE=2.74827\times10^3$}}]{\label{fig:f10} \includegraphics[scale=.183]{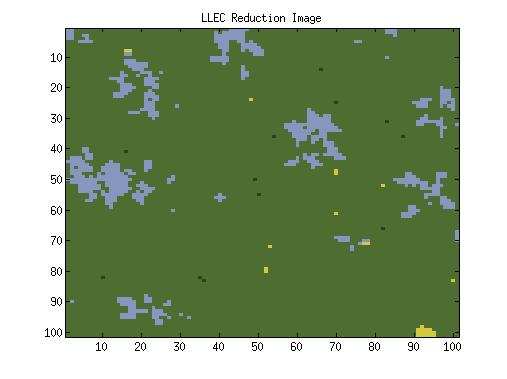}}  
  \subfloat[{\tiny{$S=8$, $DE=2.0479\times10^3$}}]{\label{fig:f11} \includegraphics[scale=.183]{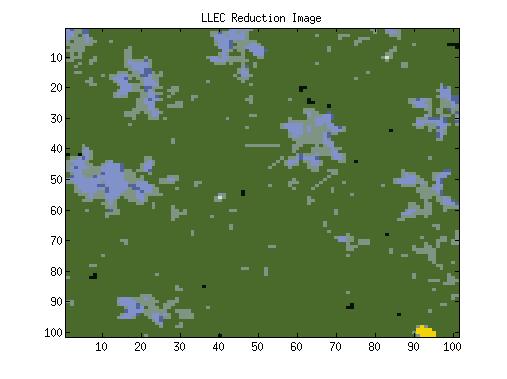}}
  \subfloat[{\tiny{$S=21$, $DE=650.323$}}]{\label{fig:f12} \includegraphics[scale=.183]{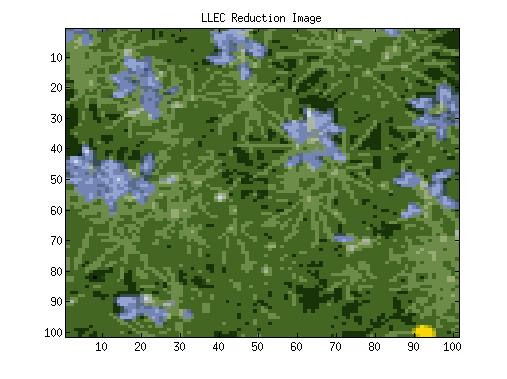}}\\
   \vspace{-.3cm}
  \subfloat[{\tiny{Original}}]{\label{fig:f13} \includegraphics[scale=.24]{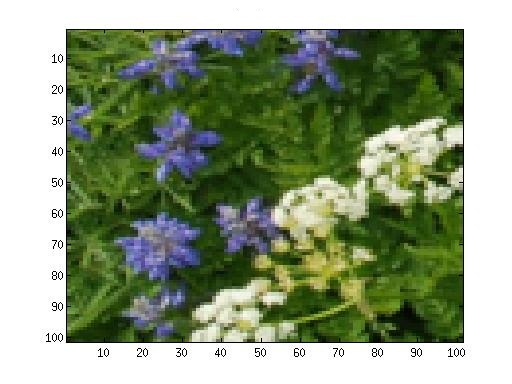}}    
  \subfloat[{\tiny{$S=5$, $DE=2.6723\times10^3$}}]{\label{fig:f14} \includegraphics[scale=.183]{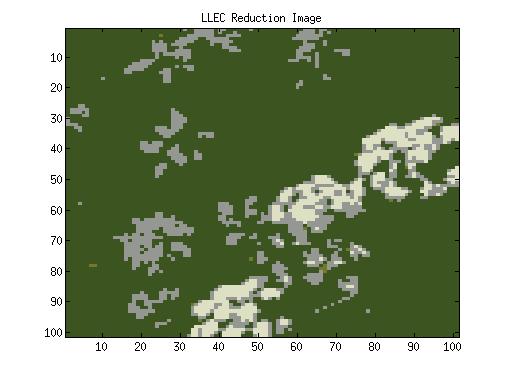}}
  \subfloat[{\tiny{$S=6$, $DE=1.6792\times10^3$}}]{\label{fig:f15} \includegraphics[scale=.183]{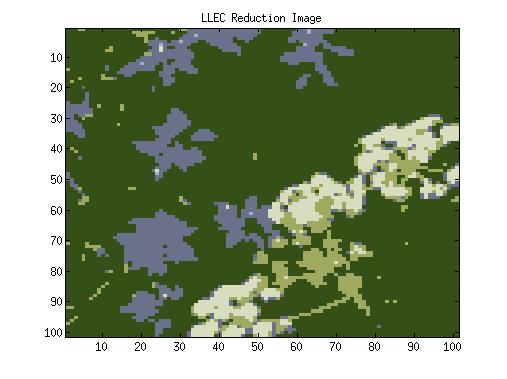}}
  \subfloat[{\tiny{$S=17$, $DE=623.756$}}]{\label{fig:f16} \includegraphics[scale=.183]{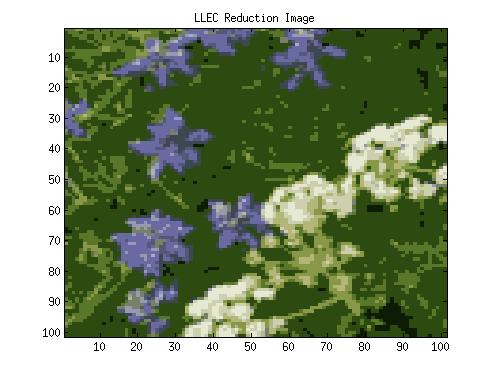}}\\
   \vspace{-.3cm}
  \subfloat[{\tiny{Original}}]{\label{fig:f17} \includegraphics[scale=.24]{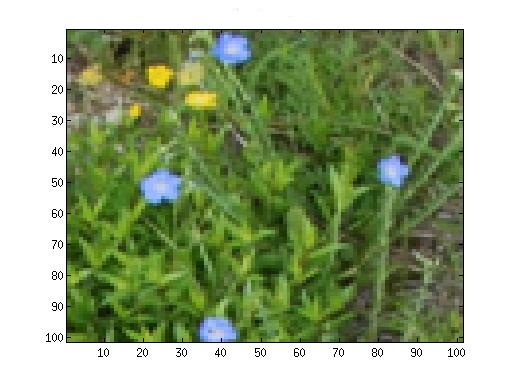}}    
  \subfloat[{\tiny{$S=4$, $DE=1.9229\times10^3$}}]{\label{fig:f18} \includegraphics[scale=.183]{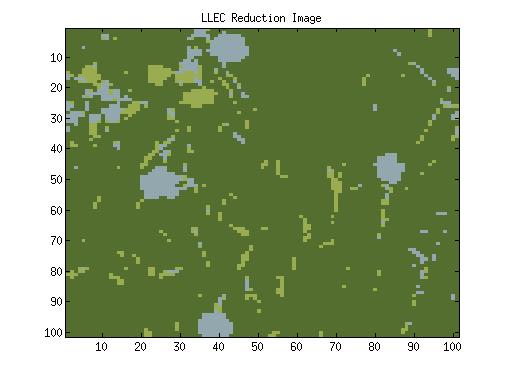}}
  \subfloat[{\tiny{$S=6$, $DE=1.0782\times10^3$}}]{\label{fig:f19} \includegraphics[scale=.183]{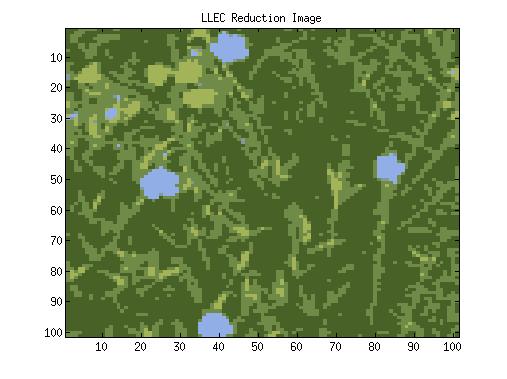}}
  \subfloat[{\tiny{$S=22$, $DE=537.293$}}]{\label{fig:f20} \includegraphics[scale=.183]{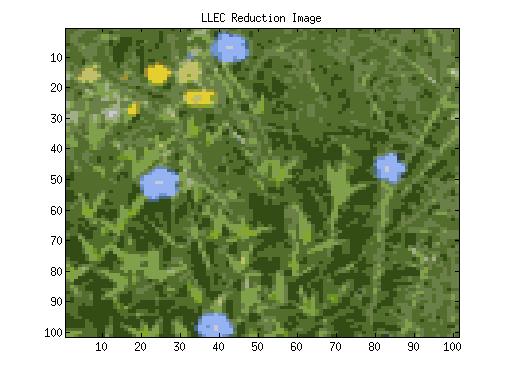}}\\
  \vspace{-.3cm}
  \caption{Reconstruction images of LLEC with variances tolerances.  Column two has tolerance $\epsilon_1=\epsilon_2=0.6$.  Column three has tolerance $\epsilon_1=\epsilon_2=0.4$.  Column four has tolerance $\epsilon_1=\epsilon_2=0.2$.  $S$ denotes the number of distinct subspaces in a reconstruction and $DE$ denotes the distortion error.}
  \label{fig:VariousFigures}
\end{figure*}

Let's first consider LLEC's ability to segment the color space of a natural image and then use this segmentation to quantize the color space.  We highlight in Figure \ref{fig:subspace} LLEC's ability to segment the subspaces in the 2-dimensional plot using accuracy tolerances of approximately 0.4 for a sample image.  In Figure \ref{fig:VariousFigures}, we observe the color quantizations obtained for this image as well as others using various accuracy tolerances.  Note that in Figure \ref{fig:VariousFigures}, each of the original images were of resolution $100 \times 100$ or less.  We see that as $\epsilon_1$ and $\epsilon_2$ decrease, the reconstructions become better representations of the original images.

\subsection{LLEC Implemented with LBG}

Let's now see how LLEC can be implemented in conjunction with another clustering algorithm on a class of large images.  These images are of a subalpine meadow near the Rocky Mountain Biological Laboratory in Gothic, Colorado provided by Dr. David Inouye of the University of Maryland.  Each image is of resolution $2592 \times 3872$ which generates 10,036,244 pixels.  We cannot implement the LLEC algorithm directly on images from this landscape data set as its implementation requires constructing a pixel by pixel matrix.  We have chosen to implement LLEC in conjuction with the Linde-Buzo-Gray algorithm \cite{Linde}, \cite{kirby} on a set of these images.  The simplicity of the LBG algorithm makes it desirable, but other clustering algorithms such as those discussed in \cite{Hartigan}, \cite{Heckbert}, \cite{kirby}, \cite{basu}, etc. could be used alternatively.  The LBG algorithm is an iterative competitive learning algorithm that, in essence, determines all points that fall within a Voronoi region around specified center vectors, calculates the mean of all points within this region, updates the center of this set to be equal to the mean, and then iterates the process until a fixed number of iterations has been met or some stopping criteria is achieved.  Proper initialization is a crucial issue for any iterative algorithm and can greatly affect the outcome.  Therefore, we have chosen four different methods to initialize the center vectors for comparison purposes.   

The first method chosen to determine the center vectors used in the iterative LBG algorithm is LLEC.  Here we use the LLEC algorithm to create a palette of colors for the data to be clustered around by identifying the natural subspaces of subimages of images within the landscape data set, namely Figures \ref{fig:f9}, \ref{fig:f13}, and \ref{fig:f17}.  The benefit of this method is that all subspaces are identified in an unsupervised manner.  

The next method involves choosing the eight three dimensional data points with components either 0 or 255 and 17 other data points sampled near the green, yellow, blue, white, and black colors as centers.  This requires supervision from the user to identify which colors seem to predominantly appear in the data set of images.

The third method involves choosing 25 random centers.  That is 25 data points, $[r,g,b]$ are chosen such that $r,g,b\in[0,255]$.  Choosing random centers in this way does not guarantee that any of the colors identified to be centers will be similar to colors that appear in the actual image, and thus, many of the centers could be potentially unused in the clustering.

The final method involves choosing 25 random centers from the data set.  That is, choose 25 columns of the data matrix $X$ randomly to be the centers that the data points are clustered around.  The benefit of this method is that this is the only approach that identifies actual points within the data set as centers.  However, not all natural subspaces may be represented as we will observe shortly. 
 \begin{figure}[t]
  \centering
\subfloat[Original]{\label{fig:L2Originalfirst} \includegraphics[scale=.48]{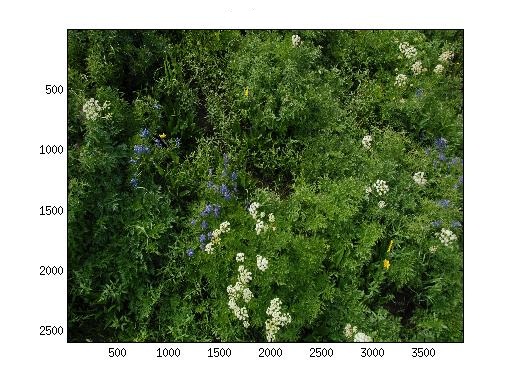}}  \\ 
\vspace{-.3cm}
  \subfloat[LLEC]{\label{fig:L2LLECfirst} \includegraphics[scale=.48]{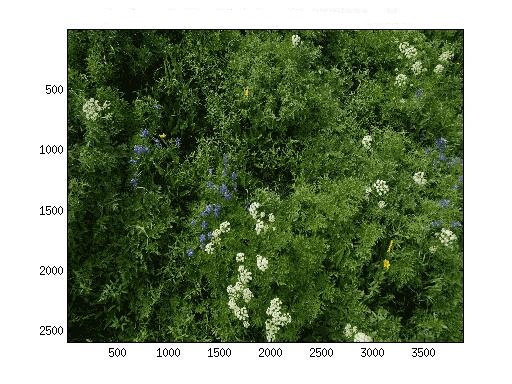}}  
  \subfloat[Identifying Centers]{\label{fig:L2Identify} \includegraphics[scale=.48]{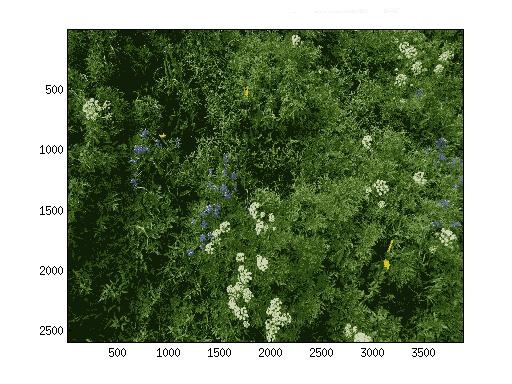}} \\
  \vspace{-.3cm}
  \subfloat[Random Centers]{\label{fig:L2Rand} \includegraphics[scale=.48]{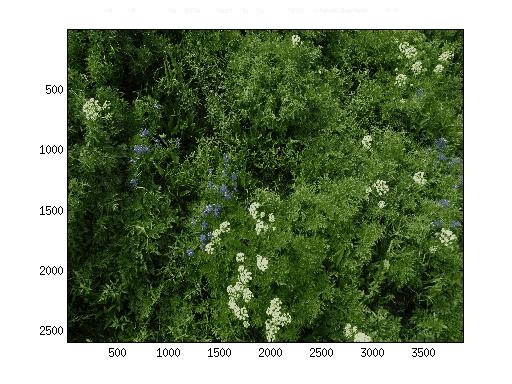}} 
  \subfloat[Random Centers from Data]{\label{fig:L2RandDatafirst} \includegraphics[scale=.48]{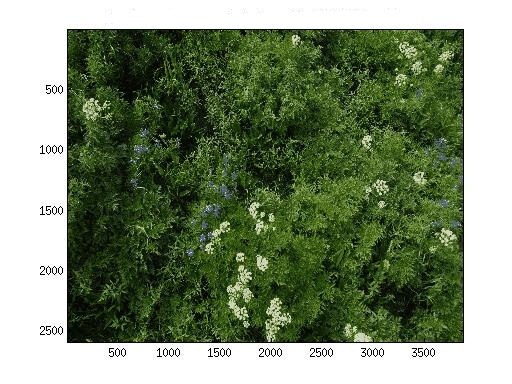}}    
  \caption{Reconstruction images after quantizing the color space of original image with LBG using indicated method to determine the centers.  Note that the respective distortion errors of each implementation with 15 iterations are:  140.0250, 342.6351, 219.0756, and 146.7013.} 
  \label{fig:LBGWithMethodsL2}
\end{figure} 

Figure \ref{fig:LBGWithMethodsL2} reveals the performance of each method on one sample image from the landscape data set.  In several implementations on various images within the landscape data set, we have observed similar results.  It appears that all methods for determining the centers result in fairly accurate reconstruction images.  However, we have observed in practice that the two methods of using the LLEC algorithm to determine centers and identifying random centers within the data set tend to result in the lowest distortion errors as calculated by 
 	$$D(X,J)=\frac{1}{p}\sum_{j \in J}\sum_{\textbf{x} \in S_j} \| \textbf{x}-\textbf{c}_j \|^2$$
where X is a data set consisting of p points with regard to a set of centers labeled by indices $J$.  Note that if we let $X^*$ indicate the matrix of points each identified with the centroid of the Voronoi region that each point is assigned to, then the distortion error could also be calculated as $\|X-X^*\|_F^2$ where $\|A\|_F =\sqrt{\sum_{i=1}^m\sum_{j=1}^n |a_{ij}|^2}$ is the Frobenius norm.  

\begin{figure*}[t]
  \centering
  \subfloat[1st Original]{\label{fig:L1Original} \includegraphics[scale=.315]{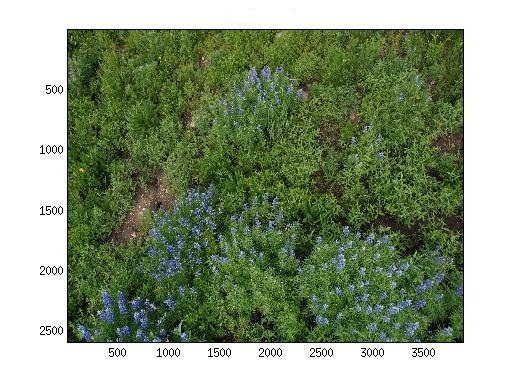}}
  \subfloat[LLEC]{\label{fig:L1LLEC} \includegraphics[scale=.315]{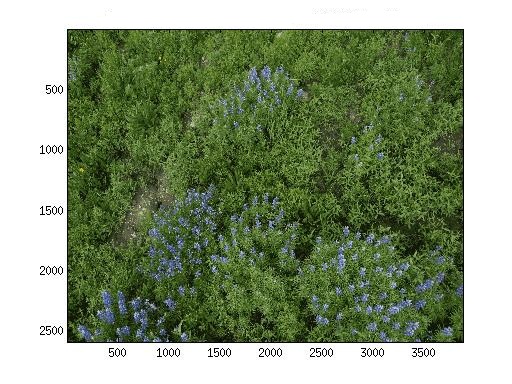}}
  \subfloat[Random from Data]{\label{fig:L1RandData} \includegraphics[scale=.315]{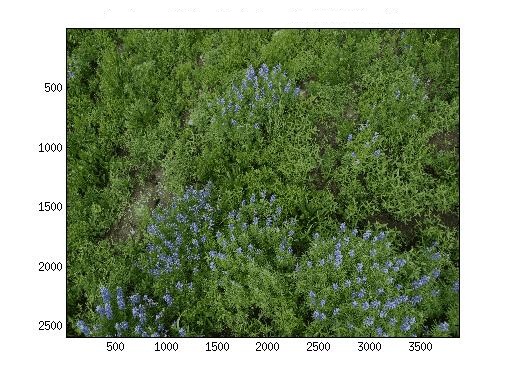}}\\   
    \subfloat[2nd Original]{\label{fig:L2Original} \includegraphics[scale=.315]{Landscape2Matlab.jpg}} 
  \subfloat[LLEC]{\label{fig:L2LLEC} \includegraphics[scale=.315]{L2LBGwithLLEC.jpg}}  
  \subfloat[Random from Data]{\label{fig:L2RandData} \includegraphics[scale=.315]{L2LBGwithRandCentersfromData2.jpg}}\\  
    \subfloat[3rd Original]{\label{fig:L3Original} \includegraphics[scale=.315]{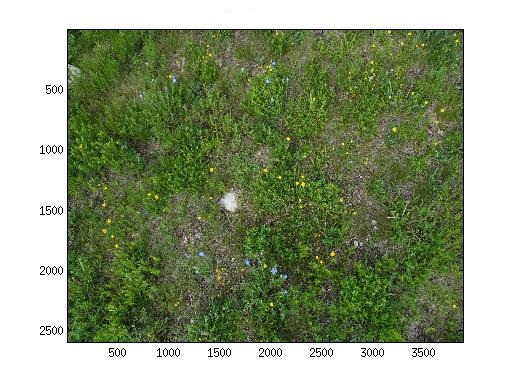}}   
  \subfloat[LLEC]{\label{fig:L3LLEC} \includegraphics[scale=.315]{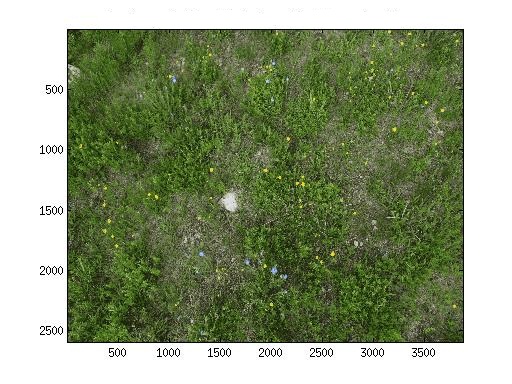}} 
  \subfloat[Random from Data]{\label{fig:L3RandData} \includegraphics[scale=.315]{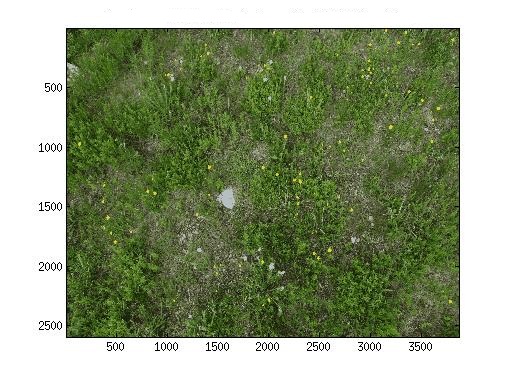}}  
  \caption{Quantizing the color space of the original image with LBG using indicated method to determine the centers.  Note that the respective distortion errors of these two implementations with 15 iterations are:  (1st Original) 210.3490 and 210.6900, (2nd Original) 140.0250 and 146.7013, (3rd Original) 172.5580 and 170.7743.} 
  \label{fig:LBGWithMethods}
\end{figure*}

We notice in Figure \ref{fig:LBGWithMethods} however that LLEC gives a better reconstruction visually.  Observe that the reconstruction obtained by identifying centers as random points within the data often does not capture all subspaces within the image.  In particular, the yellow flowers in the second example and the blue flowers in the third example do not appear in the reconstruction.  Thus, it appears that LLEC used in conjunction with LBG is able to reconstruct the image with minimal error and the most accurate representation visually.

\section{Conclusion}
In this work, we have presented a novel algorithm, LLEC, to cluster and segment the color space of natural imagery.  Within this algorithm, is a method to reconnect artificially disconnected components (resulting from a choice of $k$ nearest neighbors) as well as a technique for one dimensional subspace segmentation that can be extended to multi-dimensional segmentation which is robust in the presence of noise.  We have seen that LLEC does an excellent job of quantizing the color space of imagery with the only input parameters directly related to the accuracy of the quantization.  However, LLEC does have some limitations.  As already mentioned, LLEC is computationally intensive.  Also, in the formulation of the LLE algorithm, it is required to create matrices of size $p \times p$, where $p$ is the number of pixels in the image.  For large images, this may require a prohibitively large amount of memory.  Thus, LLE and LLEC, in turn, perform well on small sized images when being implemented in this manner.  However, if techniques such as the sampling methods discussed in \cite{Lee}, \cite{Silva} or the stitching method as discussed in \cite{bachmann} are implemented, this limitation may be alleviated.  Even with these limitations, we see that LLEC is useful in identifying the natural subspaces within an image.

\appendix


\section{Locally Linear Embedding} \label{app:lle}

\subsection{Nearest Neighbor Search} \label{sec:neighbors}

The first step in implementing the LLE algorithm is to determine the neighbors associated to each of the high dimensional data points.  Determining the nearest neighbors of a specific data point involves finding those data points that are the most similar.  One way to measure similarity is to use a Euclidean distance metric (however, other metrics may also be used).  The most straightforward way to perform this task is to determine a fixed number of nearest neighbors by considering those data points with the smallest distance determined by the metric being used.  Alternatively, nearest neighbors can be identified by classifying as neighbors all data points that fall within a ball of fixed radius about each point.  Therefore, the number of nearest neighbors could differ for each data point.  In this paper, the nearest neighbors of each point was found by determining a fixed number, $k$, of data points with the smallest non-zero Euclidean distance from the original point.  

Determining the fixed number of neighbors, $k$, is the only free parameter in the LLE algorithm.  There is some art in choosing an appropriate value for $k$.  If the manifold is well-sampled, then each data point and its nearest neighbors lie approximately on a locally linear piece of the manifold.  However, if the number of nearest neighbors, $k$, is chosen to be too large, the region may no longer be local and might include points which are geodesically far away.  The span of a set of $k$ points is a linear space of dimension at most $k-1$.  Therefore, the dimension of the target vector space, $d$, should be chosen to be strictly less than the number of nearest neighbors.  However, choosing $k$ to be too small may be problematic as the eigenvector problem to determine the embedding vectors becomes singular.  Note that \cite{saul} did not give guidance on how to choose an appropriate number of nearest neighbors.  However, \cite{kouropteva} gives an hierarchical approach to automatically select an optimal parameter value which has been shown to be quite precise.    
 
In the situation where the original dimension of the data, $D$, is fairly low, it is often necessary to choose the number of neighbors, $k$, to be greater than this dimension to avoid the eigenvector problem becoming singular.  If $k>D$, then the set of nearest neighbors, $N_i$, of data point $x_i$ is no longer linearly independent, and thus, there is not a unique solution for determining the reconstruction weights.  In this case, the covariance matrix $C$ defined below becomes singular or nearly singular.  A regularization must be implemented in order to suspend this breaking down of the algorithm.  One such regularization would be to add a small multiple of the identity to the covariance matrix which, in turn, corrects the sum of the squares of the weights so that the weights favor a uniform distribution. The optimization problem then finds the set of weights that come closest to the point representing uniform distribution of magnitude for each of the weights \cite{saul}.  In this paper, the regularization that is used is
$$C \leftarrow C + I*tol*Tr(C)$$
where $tol$ is a tolerance that is sufficiently small, usually 0.001, and $Tr(C)$ denotes the trace of $C$.  This regularizer is sufficent to make the covariance matrix well-conditioned allowing one to determine a unique solution to the optimization problem to determine the weights.

\subsection{Least Squares Problem to Find Weights}\label{sec:weights}

The second step of the LLE algorithm is to determine the weights used to associate each point with its nearest neighbors. This can be done by minimizing the distance between a point and a linear combination of all of its nearest neighbors where the coefficients of this linear combination are defined by the weights.  Let $N_i$ be the set of neighbors associated to a single point $\textbf{x}_i$, let $p$ be the number of data points being considered, and let $D$ denote the dimension of the ambient space of the data.  Our goal then is to determine the weights, $w_{ij}$, associated to each point, $\textbf{x}_i$, and each of its nearest neighbors, $\textbf{x}_j \in N_i$.  Note that the weight, $w_{ij}$, between two points that are not nearest neighbors is defined to be $0$.  Thus, a data point can only be reconstructed from points determined to be its nearest neighbors.  Now, the weights are determined by minimizing differences between each point and a linear combination of the nearest neighbors to the point.  Let $W$ denote the matrix of weights with entries $w_{ij}$. The cost function of the reconstruction error is then given by
$$\epsilon(W)= \displaystyle \sum_{i=1}^p \| \textbf{x}_i-\displaystyle\sum_{j\in{N_i}}w_{ij}\textbf{x}_j \|^2. $$
For each $i$, a constraint, $\displaystyle\sum_{j\in{N_i}}w_{ij}=1$, is implemented to ensure that these weights are invariant to translations. Note that the form of the errors ensures the weights are also invariant to rescalings and rotations. Using the sum-to-one constraint, the constraint that $w_{ij}=0$ if $x_j$ is not in the set $N_i$, and a little linear algebra, we see that
\bea
\epsilon(W)&=& \displaystyle \sum_{i=1}^p \| \textbf{x}_i 	  -\displaystyle\sum_{j\in{N_i}}w_{ij}\textbf{x}_j \|^2 
\nonumber\\
 					 &=&\displaystyle \sum_{i=1}^p \left( \displaystyle \sum_{j, k \in{N_i}}w_{ij}w_{ik}c_{jk}^i \right)
\nonumber 
\eea
where $c_{jk}^i=(\textbf{x}_i-\textbf{x}_j)^T(\textbf{x}_i-\textbf{x}_k)$.\\
Now, we want to minimize these errors using the constraint $\displaystyle\sum_{j\in{N_i}}w_{ij}=1.$  This can be done using Lagrange Multipliers.  Fixing $i$, we have  
				$$\min \displaystyle \sum_{j,k\in{N_i}}w_{j}w_{k}c_{jk}-\lambda\left( \displaystyle \sum_{j\in{N_i}}w_{j}-1\right)$$ 
This optimization problem can be solved by finding the critical values of this cost function which results in solving the following system of equations
		$$\begin{cases}
					 \displaystyle \sum_{k\in{N_i}}\tilde{w_k}c_{mk}=1\\
           \displaystyle \sum_{k\in{N_i}}\tilde{w_{k}}=1
    	\end{cases}$$
which yields 
    $$C\tilde{\textbf{w}}=\textbf{e}$$   
where $C$ corresponds to the covariance matrix determined by $c_{jk}^i=(x_i-x_j)^T(x_i-x_k)$, $\tilde{\textbf{w}}$ is the column vector of weights associated to a single point, and $\textbf{e}$ is the vector of all ones.  Thus in order to find the reconstruction weights, it is only necessary to solve
 $$\tilde{\textbf{w}}=C^{-1}\textbf{e}$$
where the weights are rescaled so that they sum to one.  Thus, we have derived the least squares problem to determine the weights that reconstruct the high dimensional data points of dimension $D$ to the lower dimension embedding data points of dimension $d$.  We can form a weight matrix, $W$, where each row, $i$, corresponds to the weights between the point $\textbf{x}_i$ and every other point.  Note that $W$ is extremely sparse as the weight between any two points that are not nearest neighbors is defined to be zero.

\subsection{Eigenvector Problem} \label{sec:e-vec}

The third and final step of the LLE algorithm is to determine the low dimensional embedding vectors, $\textbf{y}_i$, of dimension $d$ by using the reconstruction weights, $w_{ij}$, of the high dimensional data vectors, $\textbf{x}_i$.  The only information used in this portion of the algorithm is the geometry obtained by the weights.  A cost function for the errors between the reconstruction weights and the outputs, $\textbf{y}_i$, is minimized as follows:
 $$\phi(Y) =\displaystyle \sum_{i=1}^p \| \textbf{y}_i-\displaystyle \sum_{j=1}^pw_{ij}\textbf{y}_j \|^2.$$  
Here $Y$ denotes the $d \times p$ matrix of embedding vectors.  In order to find these reconstruction vectors, $\textbf{y}_i$, the following optimization problem must be solved for fixed weights, $w_{ij}$.  Using linear algebra, we can manipulate our cost function to obtain 
\bea
		\phi(Y) &=& \displaystyle \sum_{i=1}^p \| \textbf{y}_i-\displaystyle \sum_{j=1}^p w_{ij}\textbf{y}_j \|^2 \nonumber\\
						&=& \displaystyle \sum_{i,j=1}^p M_{ij} \left< \textbf{y}_i, \textbf{y}_j \right> \nonumber
\eea
where $<\ast ,\ast >$ is the standard Euclidean inner product and $M_{ij}=\delta_{ij} - w_{ij} -w_{ji} + \displaystyle \sum_{k=1}^p w_{ki}w_{kj}$.  Note that all of the $w_{ij}$ are entries of the weight matrix $W$.  Thus,
		$$M = I-W-W^T+W^TW = \left(I-W\right)^T\left(I-W\right)$$
We see that $M$ is symmetric even though $w_{ij}$ is not necessarily equal to $w_{ji}$.  In addition, $M$ is extremely sparse and is positive semi-definite.

It is straightforward to show that 
$$\phi(Y)= \displaystyle \sum_{i,j=1}^p M_{ij} \left< \textbf{y}_i, \textbf{y}_j \right> = tr(YMY^T)$$
where $Y$ corresponds to the matrix of the embedding vectors.  Our problem then becomes
		$$\min_Y tr(YMY^T)$$
		$$\text{subject to } Y Y^T=I$$
where the constraint is equivalent to saying that the embedding vectors are sphered or whitened.  Thus, they are uncorrelated, and their variances equal unity.  Note that $Y$ is orthogonal in the row space but the embedding vectors are not required to be orthogonal.  This constraint does not change the problem since the reconstruction error is invariant under rotations and rescalings.  Otherwise, letting $\textbf{y}_i=\textbf{0}$ for each $i$ would be the optimal solution.  We also use the fact that translations do not affect the cost function, so we require the outputs to be centered at the origin adding the constraint:
	$$\displaystyle \sum_{i=1}^p\textbf{\textbf{y}}_i=\textbf{0}$$
	
We will again use Lagrange multipliers to solve this problem.  Our Lagrangian becomes 
	$$L(Y,\mu)=tr(YMY^T)-\sum_{i,j=1}^d \mu_{ij} (YY^T-I)_{ij}$$
where each $\mu_{ij}$ is the Lagrange multiplier for each constraint.  Taking the derivative of this Lagrangian with respect to the matrix $Y$ and equating to zero will yield our desired solution.  We see that $2YM=2\Lambda Y$, where Lambda is the diagonal matrix of Lagrange multipliers, is the solution to our Langrangian, and it can be manipulated so that $ MY^T=Y^T\Lambda.$	
	
Thus, $Y^T$ is the matrix of eigenvectors of $M$, and $\Lambda$ is the corresponding diagonal matrix of eigenvalues.  The optimal embedding up to rotations, translations, and rescalings of the embedding space can then be found by solving this eigenvector problem.  The Rayleigh-Ritz theorem as described in \cite{horn} gives an indication of which eigenvectors actually solve the problem.  Using this, we need to obtain the bottom $d+1$ eigenvectors of the matrix, $M$, (those eigenvectors corresponding to the smallest eigenvalues in increasing order).  We will see shortly that the eigenvector corresponding to the smallest eigenvalue is the unit vector with all equal components corresponding to the mean of the data. We discard this eigenvector, leaving the second through the $d+1$ eigenvectors.  Thus, the embedding vectors that solve the LLE algorithm are these $d$ remaining eigenvectors. When discarding the bottom eigenvector, it forces each of the other eigenvectors to sum to zero by orthogonality enforcing the constraint
\bea
	\displaystyle \sum_{i=1}^p\textbf{\textbf{y}}_i=\textbf{0}
	\nonumber
\eea
which requires that the embedding vectors have zero mean.

In order for there to exist a unit eigenvector with all equal components as described above, each row of the $M$ matrix must sum to a scalar, $\lambda$.  In fact, we can show that each row has zero sum.  
Given
$$M_{ij}=\delta_{ij} -w_{ij}-w_{ji}+ \displaystyle \sum_{k=1}^p w_{ki}w_{kj}$$
we can see 
\bea 
\displaystyle \sum_{j=1}^p M_{ij}&=&\displaystyle \sum_{j=1}^p \delta_{ij} -\displaystyle \sum_{j=1}^p w_{ij}-\displaystyle \sum_{j=1}^p w_{ji}+ \displaystyle \sum_{j=1}^p \displaystyle \sum_{k=1}^p w_{ki}w_{kj}
\nonumber\\
			&=& 1 - 1 - \displaystyle \sum_{j=1}^p w_{ji}+\displaystyle \sum_{j=1}^p \displaystyle \sum_{k=1}^p w_{ki}w_{kj}
\nonumber\\
			&=& 0 - \displaystyle \sum_{j=1}^p w_{ji} + \displaystyle \sum_{k=1}^p w_{ki} 
\nonumber\\
			&=& 0
			\nonumber
\eea
Thus, there exists a unit eigenvector of all equal components corresponding to the eigenvalue zero which we may discard as described above.  See Section \ref{sec:nullspace} for a discussion on data sets that have more than one eigenvalue equal to zero.  

The third step of the LLE algorithm involves solving this eigenvector problem to determine the unique solution $Y$.  The solutions to the problem are the $d$-dimensional columns of the $Y$ matrix where each column, $j$, of $Y$ corresponds to column, $j$, of $X$, and each row of $Y$ is an eigenvector of the matrix $M$.  Note that although each weight was determined locally by reconstructing a data point by its nearest neighbors, the optimal embedding $Y$ was determined by a $p\times p$ eigensolver which is a global undertaking that uses the information from all points.  Therefore, through the LLE algorithm we were able to obtain low-dimensional embedding vectors that preserve the local topology of a high-dimensional data set by determining a global coordinate system.

\section{Initializing Subspace Segmentation} \label{app:RandPoint}

Subspace segmentation discussed in this paper is initialized by selecting a point, $\textbf{y}^*$, and then clustering around this data point.  If $\textbf{y}^*$ is selected randomly, then there is an element of randomness in the algorithm, allowing for a variety of subspaces to be determined, depending on each $\textbf{y}^*$ selected.  Here we will discuss three non-random approaches that select points reflecting the data density.

One method is to select $\textbf{y}^*$ to be the point within the data set that has the most points falling within an epsilon ball of the point.  Another method is to select $\textbf{y}^*$ to be the point whose $b^{th}$ nearest neighbor is closer to it than any other point within the data set.  Here $b \in \mathbb{Z}^+$, an arbitrary number.  While the first idea is computationally intensive, the second has the complication that points are being removed from the data set, so there may occur a moment in the algorithm where the number of data points $p<b$.  In this case, we adjust $b$ to be a number less than the number of data points. For instance, $b=\lceil \frac{p}{2} \rceil$, where $\lceil \ast \rceil$ is the ceiling function, works well in practice. 

An alternate approach for selecting the point $\textbf{y}^*$ searches for the point that falls on a subspace which most reflects a linear structure.  This can be implemented as follows.  Determine an epsilon ball around each point.  Compute the singular value decomposition of a matrix formed by the points in each of these epsilon balls in order to find the singular vectors and singular values.  Choose $\textbf{y}^*$ to be the random point with the smallest ratio of singular values $\frac{\sigma_2}{\sigma_1}$ as this reflects the subspace with the most linear structure.  While this approach chooses points falling along structures most easily identified as `linear', one downside is that it is computationally intensive.  It typically produces reconstructions with a smaller distortion error than all of the other methods described above, but it does so by identifying more subspaces.  

The various methods for identifying $\textbf{y}^*$ have features that make each of them attractive, depending on the user's desired result.  In this paper, we have chosen to follow the method that chooses $\textbf{y}^*$ as a point in a dense region of the data by finding the $b^{th}$ nearest neighbor with the smallest distance. This is the least computationally intensive and produces reconstructions with a relatively small distortion error using relatively few subspaces to reconstruct.  Here we have selected $b=50$.




\begin{thebibliography}{33}
\expandafter\ifx\csname natexlab\endcsname\relax\def\natexlab#1{#1}\fi
\providecommand{\bibinfo}[2]{#2}
\ifx\xfnm\relax \def\xfnm[#1]{\unskip,\space#1}\fi
\bibitem[{Weinberger and Saul(2004)}]{Weinberger}
\bibinfo{author}{K.~Weinberger}, \bibinfo{author}{L.~Saul}, in:
  \bibinfo{booktitle}{Computer Vision and Pattern Recognition, 2004. CVPR 2004.
  Proceedings of the 2004 IEEE Computer Society Conference on},
  volume~\bibinfo{volume}{2}, pp. \bibinfo{pages}{II--988 -- II--995 Vol.2}.
\bibitem[{Tenenbaum et~al.(2000)Tenenbaum, Silva, and Langford}]{Tenenbaum}
\bibinfo{author}{J.~B. Tenenbaum}, \bibinfo{author}{V.~Silva},
  \bibinfo{author}{J.~C. Langford}, \bibinfo{journal}{Science}
  \bibinfo{volume}{290} (\bibinfo{year}{2000}) \bibinfo{pages}{2319--2323}.
\bibitem[{Saul et~al.(2003)Saul, Roweis, and Singer}]{saul}
\bibinfo{author}{L.~K. Saul}, \bibinfo{author}{S.~T. Roweis},
  \bibinfo{author}{Y.~Singer}, \bibinfo{journal}{Journal of Machine Learning
  Research} \bibinfo{volume}{4} (\bibinfo{year}{2003})
  \bibinfo{pages}{119--155}.
\bibitem[{Kohonen(1997)}]{Kohonen}
\bibinfo{author}{T.~Kohonen}, \bibinfo{title}{Self-organizing maps},
  volume~\bibinfo{volume}{30} of \textit{\bibinfo{series}{Springer Series in
  Information Sciences}}, \bibinfo{publisher}{Springer-Verlag},
  \bibinfo{address}{Berlin}, \bibinfo{edition}{second} edition,
  \bibinfo{year}{1997}.
\bibitem[{Donoho and Grimes(2003)}]{Donoho}
\bibinfo{author}{D.~L. Donoho}, \bibinfo{author}{C.~Grimes},
  \bibinfo{title}{Hessian eigenmaps: New locally linear embedding techniques
  for high-dimensional data}, \bibinfo{year}{2003}.
\bibitem[{Huang and Mumford(????)}]{Huang}
\bibinfo{author}{J.~Huang}, \bibinfo{author}{D.~Mumford}, pp.
  \bibinfo{pages}{541--547}.
\bibitem[{Jolliffe(2002)}]{Jolliffe}
\bibinfo{author}{I.~T. Jolliffe}, \bibinfo{title}{Principal component
  analysis}, Springer Series in Statistics,
  \bibinfo{publisher}{Springer-Verlag}, \bibinfo{address}{New York},
  \bibinfo{edition}{second} edition, \bibinfo{year}{2002}.
\bibitem[{Cox and Cox(1994)}]{Cox}
\bibinfo{author}{T.~F. Cox}, \bibinfo{author}{M.~A.~A. Cox},
  \bibinfo{title}{Multidimensional scaling}, volume~\bibinfo{volume}{59} of
  \textit{\bibinfo{series}{Monographs on Statistics and Applied Probability}},
  \bibinfo{publisher}{Chapman \& Hall}, \bibinfo{address}{London},
  \bibinfo{year}{1994}. \bibinfo{note}{With 1 IBM-PC floppy disk (3.5 inch,
  HD)}.
\bibitem[{Belkin and Niyogi(2003)}]{belkin}
\bibinfo{author}{M.~Belkin}, \bibinfo{author}{P.~Niyogi},
  \bibinfo{journal}{Neural Computation} \bibinfo{volume}{15}
  (\bibinfo{year}{2003}) \bibinfo{pages}{1373--1396}.
\bibitem[{Roux et~al.(2008)Roux, Lamblin, Bengio, Joliveau, and Kégl}]{roux}
\bibinfo{author}{N.~L. Roux}, \bibinfo{author}{P.~Lamblin},
  \bibinfo{author}{Y.~Bengio}, \bibinfo{author}{M.~Joliveau},
  \bibinfo{author}{B.~Kégl}, \bibinfo{title}{Learning the 2-d topology of
  images}, \bibinfo{year}{2008}.
\bibitem[{Bachmann et~al.(2005)Bachmann, Ainsworth, and Fusina}]{bachmann}
\bibinfo{author}{C.~Bachmann}, \bibinfo{author}{T.~Ainsworth},
  \bibinfo{author}{R.~Fusina}, \bibinfo{journal}{Geoscience and Remote Sensing,
  IEEE Transactions on} \bibinfo{volume}{43} (\bibinfo{year}{2005})
  \bibinfo{pages}{441 -- 454}.
\bibitem[{Han and Goodenough(2005)}]{Han}
\bibinfo{author}{T.~Han}, \bibinfo{author}{D.~Goodenough}, in:
  \bibinfo{booktitle}{Geoscience and Remote Sensing Symposium, 2005. IGARSS
  '05. Proceedings. 2005 IEEE International}, volume~\bibinfo{volume}{2}, pp.
  \bibinfo{pages}{1237 -- 1240}.
\bibitem[{Chen et~al.(2005)Chen, Crawford, and Ghosh}]{chen}
\bibinfo{author}{Y.~Chen}, \bibinfo{author}{M.~Crawford},
  \bibinfo{author}{J.~Ghosh}, in: \bibinfo{booktitle}{Geoscience and Remote
  Sensing Symposium, 2005. IGARSS '05. Proceedings. 2005 IEEE International},
  volume~\bibinfo{volume}{6}, pp. \bibinfo{pages}{4311 -- 4314}.
\bibitem[{Westland and Ripamonti(2004)}]{Westland}
\bibinfo{author}{S.~Westland}, \bibinfo{author}{C.~Ripamonti},
  \bibinfo{title}{Computational colour science using MATLAB},
  \bibinfo{publisher}{J. Wiley}, \bibinfo{year}{2004}.
\bibitem[{Lee et~al.(2003)Lee, Pedersen, and Mumford}]{Lee}
\bibinfo{author}{A.~B. Lee}, \bibinfo{author}{K.~S. Pedersen},
  \bibinfo{author}{D.~Mumford}, \bibinfo{journal}{International Journal of
  Computer Vision} \bibinfo{volume}{54} (\bibinfo{year}{2003})
  \bibinfo{pages}{83--103}.
\bibitem[{Cheng and kuei Yang~B(2001)}]{Cheng}
\bibinfo{author}{S.-C. Cheng}, \bibinfo{author}{C.~kuei Yang~B},
  \bibinfo{title}{A fast and novel technique for color quantization using
  reduction of color space dimensionality}, \bibinfo{year}{2001}.
\bibitem[{Papamarkos et~al.(2002)Papamarkos, Atsalakis, and
  Strouthopoulos}]{Papamarkos}
\bibinfo{author}{N.~Papamarkos}, \bibinfo{author}{A.~Atsalakis},
  \bibinfo{author}{C.~Strouthopoulos}, \bibinfo{journal}{Systems, Man, and
  Cybernetics, Part B: Cybernetics, IEEE Transactions on} \bibinfo{volume}{32}
  (\bibinfo{year}{2002}) \bibinfo{pages}{44 --56}.
\bibitem[{Cheng et~al.(2001)Cheng, Jiang, Sun, and Wang}]{HDCheng}
\bibinfo{author}{H.~D. Cheng}, \bibinfo{author}{X.~H. Jiang},
  \bibinfo{author}{Y.~Sun}, \bibinfo{author}{J.~L. Wang},
  \bibinfo{journal}{Pattern Recognition} \bibinfo{volume}{34}
  (\bibinfo{year}{2001}) \bibinfo{pages}{2259--2281}.
\bibitem[{Heckbert and Heckbert(1982)}]{Heckbert}
\bibinfo{author}{P.~S. Heckbert}, \bibinfo{author}{P.~S. Heckbert},
  \bibinfo{journal}{Computer Graphics} \bibinfo{volume}{16}
  (\bibinfo{year}{1982}) \bibinfo{pages}{297--307}.
\bibitem[{Orchard and Bouman(1991)}]{Orchard}
\bibinfo{author}{M.~Orchard}, \bibinfo{author}{C.~Bouman},
  \bibinfo{journal}{Signal Processing, IEEE Transactions on}
  \bibinfo{volume}{39} (\bibinfo{year}{1991}) \bibinfo{pages}{2677 --2690}.
\bibitem[{Linde et~al.(1980)Linde, Buzo, and Gray}]{Linde}
\bibinfo{author}{Y.~Linde}, \bibinfo{author}{A.~Buzo},
  \bibinfo{author}{R.~Gray}, \bibinfo{journal}{Communications, IEEE
  Transactions on} \bibinfo{volume}{28} (\bibinfo{year}{1980})
  \bibinfo{pages}{84 -- 95}.
\bibitem[{Deng and Manjunath(2001)}]{Deng}
\bibinfo{author}{Y.~Deng}, \bibinfo{author}{B.~Manjunath},
  \bibinfo{journal}{Pattern Analysis and Machine Intelligence, IEEE
  Transactions on} \bibinfo{volume}{23} (\bibinfo{year}{2001})
  \bibinfo{pages}{800 --810}.
\bibitem[{Velho et~al.(1997)Velho, Gomes, Vinicius, and Sobreiro}]{Velho}
\bibinfo{author}{L.~Velho}, \bibinfo{author}{J.~Gomes},
  \bibinfo{author}{M.~Vinicius}, \bibinfo{author}{R.~Sobreiro}, in:
  \bibinfo{booktitle}{in Proc. Tenth Brazilian Symp. Comput. Graph. Image
  Process}, \bibinfo{publisher}{IEEE Computer Society}, \bibinfo{year}{1997},
  pp. \bibinfo{pages}{203--210}.
\bibitem[{Ng et~al.(2001)Ng, Jordan, and Weiss}]{ng}
\bibinfo{author}{A.~Y. Ng}, \bibinfo{author}{M.~I. Jordan},
  \bibinfo{author}{Y.~Weiss}, in: \bibinfo{booktitle}{ADVANCES IN NEURAL
  INFORMATION PROCESSING SYSTEMS}, \bibinfo{publisher}{MIT Press},
  \bibinfo{year}{2001}, pp. \bibinfo{pages}{849--856}.
\bibitem[{Kirby(2001)}]{kirby}
\bibinfo{author}{M.~Kirby}, \bibinfo{title}{Geometric data analysis},
  \bibinfo{publisher}{Wiley-Interscience [John Wiley \& Sons]},
  \bibinfo{address}{New York}, \bibinfo{year}{2001}. \bibinfo{note}{An
  empirical approach to dimensionality reduction and the study of patterns}.
\bibitem[{Chung(1997)}]{Chung}
\bibinfo{author}{F.~R.~K. Chung}, \bibinfo{title}{Spectral graph theory},
  volume~\bibinfo{volume}{92} of \textit{\bibinfo{series}{CBMS Regional
  Conference Series in Mathematics}}, \bibinfo{publisher}{Published for the
  Conference Board of the Mathematical Sciences, Washington, DC},
  \bibinfo{year}{1997}.
\bibitem[{Basu et~al.(2008)Basu, Davidson, and Wagstaff}]{basu}
\bibinfo{author}{S.~Basu}, \bibinfo{author}{I.~Davidson},
  \bibinfo{author}{K.~Wagstaff}, \bibinfo{title}{Constrained Clustering:
  Advances in Algorithms, Theory, and Applications},
  \bibinfo{publisher}{Chapman \& Hall/CRC}, \bibinfo{edition}{1} edition,
  \bibinfo{year}{2008}.
\bibitem[{Shi and Malik(1997)}]{shi}
\bibinfo{author}{J.~Shi}, \bibinfo{author}{J.~Malik}, \bibinfo{journal}{IEEE
  Transactions on Pattern Analysis and Machine Intelligence}
  \bibinfo{volume}{22} (\bibinfo{year}{1997}) \bibinfo{pages}{888--905}.
\bibitem[{Vidal et~al.(2005)Vidal, Ma, and Sastry}]{Vidal}
\bibinfo{author}{R.~Vidal}, \bibinfo{author}{Y.~Ma},
  \bibinfo{author}{S.~Sastry}, \bibinfo{journal}{Pattern Analysis and Machine
  Intelligence, IEEE Transactions on} \bibinfo{volume}{27}
  (\bibinfo{year}{2005}) \bibinfo{pages}{1945 --1959}.
\bibitem[{Hartigan(1975)}]{Hartigan}
\bibinfo{author}{J.~A. Hartigan}, \bibinfo{title}{Clustering algorithms},
  \bibinfo{publisher}{John Wiley \& Sons, New York-London-Sydney},
  \bibinfo{year}{1975}. \bibinfo{note}{Wiley Series in Probability and
  Mathematical Statistics}.
\bibitem[{de~Silva and Carlsson(2004)}]{Silva}
\bibinfo{author}{V.~de~Silva}, \bibinfo{author}{G.~Carlsson},
  \bibinfo{journal}{IEEE Symposium on Point-based Graphic}
  (\bibinfo{year}{2004}) \bibinfo{pages}{157--166}.
\bibitem[{Kouropteva et~al.(2002)Kouropteva, Okun, and
  Pietikäinen}]{kouropteva}
\bibinfo{author}{O.~Kouropteva}, \bibinfo{author}{O.~Okun},
  \bibinfo{author}{M.~Pietikäinen}, in: \bibinfo{booktitle}{1 st International
  Conference on Fuzzy Systems and}, pp. \bibinfo{pages}{359--363}.
\bibitem[{Horn and Johnson(1990)}]{horn}
\bibinfo{author}{R.~A. Horn}, \bibinfo{author}{C.~R. Johnson},
  \bibinfo{title}{Matrix analysis}, \bibinfo{publisher}{Cambridge University
  Press}, \bibinfo{address}{Cambridge}, \bibinfo{year}{1990}.
  \bibinfo{note}{Corrected reprint of the 1985 original}.

\end{thebibliography}







\end{document}